# Трехстадийная модель равновесного распределения транспортных потоков

*А.В. Гасников (МФТИ),*
*Ю.В. Дорн (МФТИ),*
*Ю.Е. Нестеров (ВШЭ),*
*С.В. Шпирко (МФТИ)*

В работе предпринята попытка объединить в одну общую модель, сводящуюся к решению задачи негладкой выпуклой оптимизации: модель расчета матрицы корреспонденций (энтропийная модель), модель равновесного расщепления потоков (по способу передвижений) и модель равновесного распределения потоков по путям (модель стабильной динамики, также называемую в литературе моделью Нестерова–деПальма, 1998).

Данная работа представляет собой переработанную и осовремененную (в контексте недавно полученных численных результатов) версию работы:

*Гасников А.В., Дорн Ю.В., Нестеров Ю.Е, Шпирко С.В.* О трехстадийной версии модели стационарной динамики транспортных потоков // Математическое моделирование. 2014. Т. 26:6. С. 34–70.

**Ключевые слова:** равновесное распределение транспортных потоков, расщепление потоков, матрица корреспонденций, модель стабильной динамики, теорема о минимаксе, принцип Вардропа, равновесие Нэша.
Библиография: 65 названий, 3 рисунка.

Работа А.В. Гасникова выполнена при поддержке гранта РФФИ 15-31-70001 мол_а_мос. Работа Ю.В. Дорна выполнена при поддержке гранта РФФИ 15-31-20571-мол_а_вед. Работа Ю.Е. Нестерова и С.В. Шпирко выполнена при поддержке гранта РФФИ 14-01-00722-а.

## Оглавление

# 1. Введение

Одной из основных задач последнего времени, остро стоящих в Москве и ряде других крупных городов России (Санкт-Петербург, Пермь, Владивосток, Иркутск, Калининград и др.) является разработка транспортной модели города, позволяющей решать задачи долгосрочного планирования (развития) транспортной инфраструктуры города [1]. В частности, ожидается, что разработка такой модели поможет ответить на вопросы: какой из проектов дорожного строительства оптимален, где пропуская способность дороги недостаточна, как изменится транспортная ситуация, если построить в этом месте торговый центр (жилой район, стадион), как правильно определять маршруты и расписание движения общественного транспорта, какой эффект даст выделение полос для общественного транспорта и т.п. Уже имеется программное обеспечение, позволяющее частично решать указанные выше задачи. Однако имеется много вопросов к тому, какие модели и алгоритмы используются в большинстве программных продуктах. Например, не очевидным элементом почти всех этих продуктов является использование в качестве одного из блоков модели равновесного распределения транспортных потоков Бэкмана (1955) [2–6]. Эта во многом хорошая модель, тем не менее, имеет довольно много недостатков (см. [7]). Например, калибровка такой модели требует знания функций затрат на ребрах графа транспортной сети (эти функции связывают время в пути по ребру с величиной транспортного потока по этому ребру), причем сама модель оказывается довольно чувствительной к выбору этих функций, которые в модели Бэкмана, как правило, предполагаются выпуклыми, монотонно возрастающими. Достаточно сказать, что в случае наличия платных дорог, для расчета оптимальных плат за проезд требуется вычислять, например, производные этих неизвестных функций [8, 9]. Существование таких функций в модели Бэкмана является одним из основных предположений, и, одновременно, одним из самых слабых мест. Реальные данные показывают (см. рис. 1, полученный В.А. Данилкиным по данным ЦОДД в 2012 году), что предположение о классе функций затрат не выполняется. Но даже если предположить, что такая зависимость все же существует,[1] то по-прежнему остается другая проблема: как калибровать модель Бэкмана, то есть откуда брать эти зависимости. Не получится ли переобучения у создаваемой нами модели? То есть, не получится ли так, что распоряжаясь большим произволом при калибровке по обучающей выборке (историческим данным) мы “переподгоним” модель: исторические данные за счет большого числа подкручиваемых параметров мы, действительно, можем хорошо научиться описывать, но использовать такую модель для планирования будет опасно, поскольку не будет контроля переобучения. Обычным средством борьбы с переобучением в этом месте является параметризация функций затрат, например, в классе BPR функций [1–6]. К сожалению, какого бы то ни было научного обоснования, почему именно такая параметризация используется, нам не известно.

[1] Это допущение можно оправдать, например, тем, что, как правило, мы рассматриваем равновесные конфигурации с точки зрения пользователей сети в модели Бэкмана, которым соответствует только одна из веток – верхняя (рис. 1). В модели Бэкмана пользователи сети при принятии решений оценивают время в пути в зависимости от величины потока, которая интерпретируется, как число желающих проехать по этому ребру в единицу времени. Нижняя ветка, отвечающая приблизительно линейному росту скорости с ростом потока $V \approx f/\rho_{\max}$, соответствует ситуации, когда есть узкое место, пропускная способность которого по каким-то причинам определяется не типичными характеристиками рассматриваемого ребра, а скажем, пробкой, пришедшей с впереди идущего ребра (по ходу движения). И в таких ситуациях величина потока $f$ интерпретируется не как число желающих, а как число могущих проехать. Такие ситуации просто исключаются в модели Бэкмана.

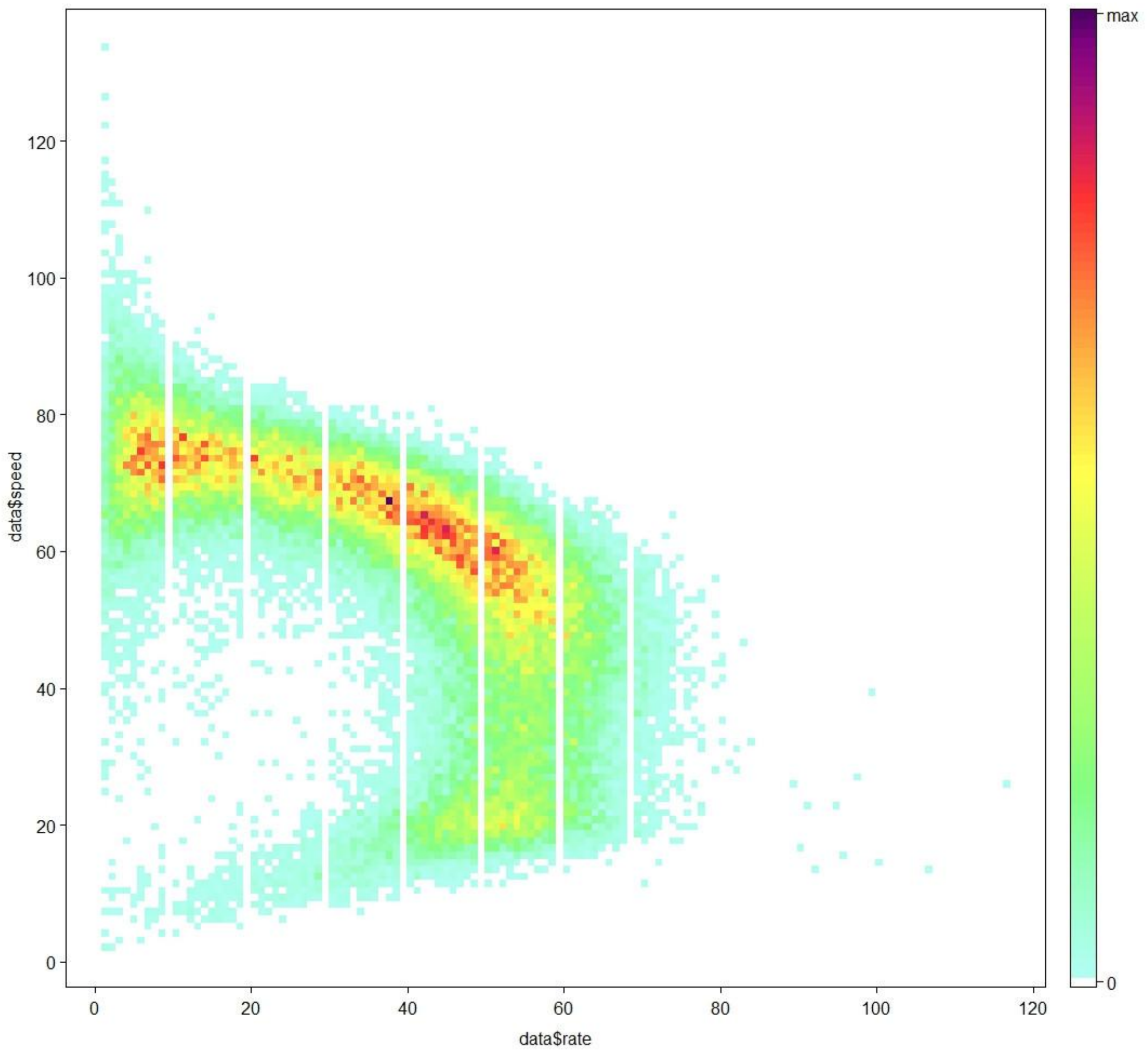

По оси абсцисс – поток по двум полосам (авт/мин), а по оси ординат – скорость (км/час)

**Рис. 1**

Другим, не очевидным элементом этих программных продуктов являются используемые вычислительные алгоритмы: контроль их робастности [10], к неточности (неполноте) данных, ошибкам округления (поскольку возникают задачи огромных размеров, то такие ошибки могут интенсивно накапливаться). Наконец, сама философия, использующаяся в таких продуктах при построении равновесной модели города [1] также вызывает много вопросов, о которых немного подробнее будет написано в п. 5. Несмотря на отмеченные выше проблемы, разработчики программного обеспечения часто находят вполне разумные инженерные компромиссы (нам известно о положительном опыте PTV VISSUM и TRANSNET), не плохо работающие практике.

Целью данной работы является предложить математическую трехстадийную транспортную модель города, в которой один из блоков (модель равновесного распределения потоков) предлагается заменить с модели Бэкмана на модель стабильной динамики [7, 11]. К сожалению, целый ряд проблем, свойственных ранее известным многостадийным моделям будет присущ и модели, предлагаемой в данной работе. Однако несколько важных недостатков, по-видимому, удалось устранить. Прежде всего, речь идет о возможности калибровки модели по реальным данным, контроле переобучения и существовании эффективного робастного вычислительного алгоритма с гарантированными (не улучшаемыми) оценками числа затраченных арифметических операций для достижения требуемой точности. Последнее обстоятельство представляется особенно важным в контексте того, как обычно используются такие модели. А именно, с помощью таких моделей просматривается множество различных сценариев. К сожалению, оптимизационные задачи вида: “где и какую дорогу стоит построить при заданных бюджетных ограничениях”, а также многие другие задачи, решаются непосредственным перебором различных вариантов, где для расчета каждого варианта потребуется запускать модель, меняя каждый раз что-то на входе. Кстати сказать, предложенный для расчета модели алгоритм позволяет также учитывать масштаб изменения входных данных. Если эти изменения не большие (точнее говоря, меняется не большое количество входных параметров), то для выполнения перерасчета по модели потребуется значительно меньше времени, чем при первом запуске.

## 2. Структура работы и предварительные сведения

Опишем вкратце структуру работы. В п. 3 описывается эволюционный способ вывода популярного на практике статического способа расчета матрицы корреспонденций (энтропийной модели). Также приводится основная идея, базирующаяся на теореме Тихонова о разделении времен [12], получения трехстадийной модели из отдельных блоков (расчет матрицы корреспонденций + равновесное расщепление потоков + равновесное распределение потоков). Точнее говоря, использование теоремы Тихонова – лишь часть идеи, которая сводит решение задачи к поиску (единственного) притягивающего положения равновесия системы в медленном времени, отвечающей за формирование корреспонденций, при подстановке в неё зависимостей времен в пути от корреспонденций (такие зависимости получаются из системы в быстром времени). Другая её часть, заключается в том, что задачи поиска (единственного) притягивающего положения равновесия системы в медленном времени и поиска зависимостей времен в пути от корреспонденций с помощью некоторых вариационных принципов, о которых говорится в пп. 3 и 4, 6–8, сводятся к задачам выпуклой оптимизации, которые можно объединить в одну общую задачу поиска седловой точки негладкой выпукло-вогнутой функции. В конце пункта указывается возможность обобщения трехстадийной модели до четырехстадийной (есть некоторые нюансы в трактовках – стоит иметь в виду, что в разных литературных источниках к таким моделям могут предъявляться немного разные требования), в которой учитываются различные типы пользователей и различные типы передвижений (моделирование идет на больших масштабах времени). В п. 4 описывается модель равновесного распределения потоков Бэкмана. В конце пункта приводится эволюционный способ интерпретации возникающего в этой модели равновесия Нэша–Вардропа. В п. 5 приводится краткий обзор многостадийных моделей, построенных на основе моделей описанных в пп. 3 и 4. В п. 6 описывается модель равновесного распределения потоков, которую мы далее будем называть моделью стабильной динамики. Модель стабильной динамики требует намного меньше данных для своей калибровки, наследует практически все основные “хорошие” свойства модели Бэкмана, и не наследует ряд недостатков. В п. 7 модель стабильной динамики выводится с помощью предельного перехода из модели Бэкмана. В п. 8 строится обобщение модели стабильной динамики на случай, когда есть несколько способов передвижения (личный транспорт и общественный; отметим при этом, что в Москве и области более 70% пользователей сети используют общественный транспорт). Таким образом, в пункте 8 в модель стабильной динамики органично встраивается модель равновесного расщепления потоков. В п. 9 энтропийная модель расчета матрицы корреспонденций из п. 3 объединяется с моделью из п. 8. В результате получается трехстадийная модель, в которой учитывается и формирование корреспонденций, и расщепление потоков, и равновесное распределение потоков по графу транспортной сети. Примечательно, что поиск равновесия в полученной трехстадийной модели (из задачи поиска седловой точки негладкой выпукло-вогнутой функции – см. п. 3) в итоге сводится к решению задачи негладкой выпуклой оптимизации, с ограниченной константой Липшица функционала, но с неконтролируемой начальной невязкой. Важно отметить, что помимо самого решения, нужно определять и часть двойственных переменных, имеющих содержательный физический смысл. В п. 10 модель п. 9 обобщается на случай поиска стохастического равновесия, что можно проинтерпретировать как ограниченную рациональность водителей или их неполную информированность. Это допущение делает модель п. 9 более приближенной к практике. С вычислительной точки зрения, сделанная модификация сводит задачу к задаче гладкой выпуклой оптимизации с немного более громоздким функционалом. Полученная задача во многом наследует все вычислительные минусы и плюсы негладкого случая. Обе задачи выпуклой оптимизации из пп. 9, 10 требуют разработки адекватных прямо-двойственных субградиентных алгоритмов решения. Заметим, что использовать методы с оракулом (сообщающим, в зависимости от своего порядка, значения функций в выбранных точках, их градиенты, и т.д.) порядка выше первого [13] не представляется возможным в виду размеров задач. В заключительном пункте 11 приводятся (с объяснениями) несколько практических рецептов по калибровке предложенных моделей.

На протяжении всей работы (и последующих работ в этом направлении) мы будем активно использовать элементы выпуклого анализа и методы численного решения задач выпуклой оптимизации ориентируясь на априорное знакомство читателя с этими дисциплинами, например, в объеме книг [13, 14]. Одним из основных инструментов работы будет теорема фон Неймана о минимаксе для выпукло-вогнутых функций [15]. Причем использоваться эта теорема будет не только для функций, заданных на произведении компактов, но и на неограниченных множествах. Надо лишь иметь гарантию, что максимумы и минимумы существуют (достигаются). Проблема сводится к существованию неподвижной точки у многозначного отображения. В этой области имеется большое количество результатов, с запасом покрывающих потребности данной работы, для обоснования возможности перемены порядка взятия максимума и минимума. В частности, в этой работе мы будем пользоваться вариантом минимаксной теоремы, называемой в зарубежной литературе “Sion’s minimax theorem” [16], в которой предполагается компактность лишь одного из множеств, отвечающих выпуклым или вогнутым переменным, и непрерывность функции. Тем не менее, далее мы будем использовать более привычное название “минимаксная теорема фон Неймана”.

Все необходимые обозначения и ссылки будут вводиться в работе по мере необходимости.

## 3. Энтропийная модель расчета матрицы корреспонденций

Приведем, во многом следуя книге [6], обоснование (ориентированное, скорее, на западный уклад жизни, поскольку предполагается, что человек за жизнь несколько раз меняет место жительство и место работы), пожалуй, одному из самых популярных способов расчета матрицы корреспонденций, имеющему более чем сорокалетнюю историю, – энтропийной модели [17–19].

Пусть в некотором городе имеется $n$ районов, $L_i > 0$ – число жителей $i$-го района, $W_j > 0$ – число работающих в $j$-м районе. При этом $N = \sum_{i=1}^{n} L_i = \sum_{j=1}^{n} W_j$ , – общее число жителей города. В последующих пунктах, под $L_i \geq 0$ будет пониматься число жителей района, выезжающих в типичный день за рассматриваемый промежуток времени из $i$-го рай-она, а под $W_j \geq 0$ – число жителей города, приезжающих на работу в $j$-й район в типичный день за рассматриваемый промежуток времени. Обычно, так введенные, $L_i$, $W_j$ рассчитываются через число жителей $i$-го района и число работающих в $j$-м районе с помощью более менее универсальных (в частности, не зависящих от $i$, $j$) коэффициентов пропорциональности. Эти величины являются входными параметрами модели, т.е. они не моделируются (во всяком случае, в рамках выбранного подхода). Для долгосрочных расчетов с разрабатываемой моделью требуется иметь прогноз изменения значений этих величин.

Обозначим через $d_{ij}(t) \geq 0$ – число жителей, живущих в $i$-м районе и работающих в $j$-м в момент времени $t$. Со временем жители могут только меняться квартирами, поэтому во все моменты времени $t \geq 0$

$$d_{ij}(t) \geq 0,\ \sum_{j=1}^{n} d_{ij}(t) \equiv L_i,\ \sum_{i=1}^{n} d_{ij}(t) \equiv W_j,\ i, j = 1, ..., n,\ 1 \ll n^2 \ll N. \quad \textbf{(A)}$$ [2]

Опишем основной стимул к обмену: работать далеко от дома плохо из-за транспортных издержек. Будем считать, что эффективной функцией затрат [6] будет $R(T) = \beta T/2$, где $T > 0$ – время в пути от дома до работы, а $\beta > 0$ – настраиваемый параметр модели (который также можно проинтерпретировать и даже оценить).

Теперь опишем саму динамику (детали см. в [6, 21]). Пусть в момент времени $t \geq 0$ $r$-й житель живет в $k$-м районе и работает в $m$-м, а $s$-й житель живет в $p$-м районе и работает в $q$-м. Тогда $\lambda_{k,m;\,p,q}(t)\Delta t + o(\Delta t)$ – есть вероятность того, что жители с номерами $r$ и $s$ ($1 \leq r < s \leq N$) "поменяются" квартирами в промежутке времени $(t, t+\Delta t)$. Вероятность обмена местами жительства зависит только от мест проживания и работы обменивающихся:

$$\lambda_{k,m;\,p,q}(t) \equiv \lambda_{k,m;\,p,q} = \lambda N^{-1} \exp\Big( \underbrace{R(T_{km}) + R(T_{pq})}_{\text{суммарные затраты до обмена}} - \underbrace{\big(R(T_{pm}) + R(T_{kq})\big)}_{\text{суммарные затраты после обмена}} \Big) > 0,$$

---

[2] Для Москвы (и других крупных мегаполисов) часто выбирают $n \sim 10^2 - 10^3$. Следовательно, корреспонденций будет $n^2 \sim 10^4 - 10^6$, и чтобы каждую из корреспонденций определить с точностью 10% (относительной точностью $\varepsilon \sim 10^{-1}$) потребуется (это оценка снизу) опросить не менее $n^2/\varepsilon^2 \sim 10^6 - 10^8$ жителей города, что не представляется возможным – это обстоятельство является одной из причин (не единственной), почему стараются снизить размерность пространства параметров, считая, что матрица корреспонденций задается не $n^2$ параметрами, а только, например, $2n$ ($L$ и $W$), которые и надо определять (см. также п. 5, в котором параметров $2n+1$). На самом деле выбирать критерием качества "относительную точность" для каждой (в том числе очень маленькой) корреспонденции не очень разумно. Более естественно восстанавливать матрицу (вектор) корреспонденций в 1-норме. Другой причиной использования описанной ниже энтропийной модели является возможность прогнозирования с её помощью того, как будет меняться матрица корреспонденций при изменении инфраструктуры города, собственно эта одна из тех задач, которые необходимо уметь решать, для получения ответов на вопросы, приведенные в п. 1.

Кроме того, важно заметить, что крупный Мегаполис, как правило, представляется в такого рода моделях вместе со всеми своими окрестными территориями. Скажем, для Москвы – это Московская область. Отметим также, что для Москвы и области $N \sim 10^7$.

где коэффициент $\lambda > 0$ характеризует интенсивность обменов. Совершенно аналогичным образом можно было рассматривать случай "обмена местами работы". Здесь стоит оговориться, что "обмены" не стоит понимать буквально – это лишь одна из возможных интерпретаций. Фактически используется, так называемое, "приближение среднего поля" [20], т.е. некое равноправие агентов (жителей) внутри фиксированной корреспонденции и их независимость.[3]

Согласно эргодической теореме для марковских цепей (в независимости от начальной конфигурации $\left\{d_{ij}(0)\right\}_{i=1,\,j=1}^{n,\,n}$) [6, 24] предельное распределение совпадает со стационарным (инвариантным), которое можно посчитать (получается проекция прямого произведение распределений Пуассона на гипергрань некоторого многогранника):

$$\exists\ c_n > 0:\ \forall\ \left\{d_{ij}\right\}_{i=1,\,j=1}^{n,\,n} \in (\mathrm{A}),\ t \ge c_n N \ln N$$

$$P\left(d_{ij}(t) = d_{ij}, i, j = 1,\ldots,n\right) \simeq Z^{-1}\prod_{i,j=1}^{n}\exp\left(-2R\left(T_{ij}\right)d_{ij}\right)\cdot\left(d_{ij}!\right)^{-1} \overset{def}{=} p\left(\left\{d_{ij}\right\}_{i=1,\,j=1}^{n,\,n}\right),$$

где "статсумма" $Z$ находится из условия нормировки получившейся "пуассоновской" вероятностной меры. Отметим, что стационарное распределение $p\left(\left\{d_{ij}\right\}_{i=1,\,j=1}^{n,\,n}\right)$ удовлетворяет условию детального равновесия:

$$\left(d_{km}+1\right)\left(d_{pq}+1\right)p\left(\left\{d_{11},\ldots,d_{km}+1,\ldots,d_{pq}+1,\ldots,d_{pm}-1,\ldots,d_{kq}-1,\ldots,d_{nn}\right\}\right)\lambda_{k,m;\,p,q} =$$

$$= d_{pm}d_{kq}\,p\left(\left\{d_{ij}\right\}_{i=1,\,j=1}^{n,\,n}\right)\lambda_{p,m;\,k,q}.$$

При $N \gg 1$ распределение $p\left(\left\{d_{ij}\right\}_{i=1,\,j=1}^{n,\,n}\right)$ экспоненциально сконцентрировано на множестве (**A**) в $\mathrm{O}\left(\sqrt{N}\right)$ окрестности наиболее вероятного значения $\left\{d_{ij}^{*}\right\}_{i=1,\,j=1}^{n,\,n}$, которое определяется, как решение задачи энтропийно-линейного программирования [6] (подробнее о таких задачах и о том, как их решать см. [17–20, 25–33] – при небольших $\beta$ успешно работает метод балансировки [29], при больших $\beta$ доминирует подход из [32, 33]):[4]

$$\ln p\left(\left\{d_{ij}\right\}_{i=1,\,j=1}^{n,\,n}\right) \sim -\sum_{i,j=1}^{n} d_{ij}\ln\left(d_{ij}/e\right) - \beta\sum_{i,j=1}^{n} d_{ij}T_{ij} \to \max_{\left\{d_{ij}\right\}_{i=1,\,j=1}^{n,\,n}\in(\mathrm{A})}. \tag{1}$$

---

[3] Отметим также, что конечной цели (получение задачи (1)) можно добиться отличными способами. Скажем, используя формализм Л.И. Розоноэра "систем обмена и распределения ресурсов" со структурной функцией энтропией (кстати, это можно не постулировать, структурная функция появляется, при работе с условием интегрируемости дифференциальной формы возможных обменов ресурсами) – см., например, [22] и цитированную там литературу. При таком подходе стохастика не нужна, и вариационный принцип (максимизации энтропии при аффинных ограничениях) получается мало чувствительным к особенностям возможных превращений в системе. Другие способы получения (1) связаны с информационно-статистическими соображениями, например, принципом максимума правдоподобия [17, 19, 23].

[4] При получении этой формулы использовалась асимптотическая формула Стирлинга [6, 17, 19], то есть предполагалось, что если $d_{ij} > 0$, то $d_{ij} \gg 1$ (если все $L_i > 0$, $W_j > 0$, то и все $d_{ij} > 0$ [6, 19]) и, как следствие, целочисленностью переменных $d_{ij}$ можно пренебречь, то есть решать не *NP*-полную задачу выпуклого целочисленного программирования, а обычную задачу выпуклой оптимизации (1). Сделанное предположение "если $d_{ij} > 0$, то $d_{ij} \gg 1$" во многом будет следовать из дальнейших рассуждений (см. также [6, 19]).

Естественно принимать решение этой задачи $\left\{d_{ij}^*\right\}_{i=1,\, j=1}^{n,\, n}$ за равновесную конфигурацию [34]. Однако имеется проблема: $T_{ij}$ – неизвестны, и зависят от $\left\{d_{ij}\right\}$. Эту проблему мы постараемся решить в дальнейшем.

Обратим внимание, что предложенный выше вывод энтропийной модели расчета матрицы корреспонденций отличается от классического [17]. В монографии А.Дж. Вильсона [17] $\beta$ интерпретируется как множитель Лагранжа к ограничению на среднее "время в пути": $\sum_{i,j=1}^{n} d_{ij}T_{ij} = C$.[5] Сделано это нами для того, чтобы контролировать знак параметра $\beta > 0$ и лучше понимать его физический смысл (нам это понадобится в дальнейшем).[6] Отметим также, что характерный временной масштаб формирования корреспонденций – годы (это не совсем так для корреспонденций типа дом–торговля, дом–отдых). В то время как характерное время установления равновесных значений $T_{ij}(d)$ недели (см. ниже). Хочется сказать, что мы здесь находимся в условиях теоремы Тихонова о разделении времен [12], точнее говоря, что можно просто подставить в (1) зависимости $T_{ij}(d)$ и решать полученную задачу. Но теорема Тихонова должна применяться для системы ОДУ, представляющей собой, в данном случае, динамику квазисредних введенной выше стохастической динамики [6, 20]:[7]

$$\frac{d}{dt}c_{ij} = \sum_{k,p=1}^{n} \lambda \exp\left(\frac{\beta}{2}\left(\left[T_{ip} + T_{kj}\right] - \left[T_{ij} + T_{kp}\right]\right)\right) c_{ip}c_{kj} -$$

$$-\sum_{k,p=1}^{n} \lambda \exp\left(\frac{\beta}{2}\left(\left[T_{ij} + T_{kp}\right] - \left[T_{ip} + T_{kj}\right]\right)\right) c_{ij}c_{kp},\ c_{ij} = d_{ij}/N;$$

$$\varepsilon \frac{d}{dt}T_{ij} = \left[\text{сложный оператор, зависящий от } c\right],\ \varepsilon \sim 10^{-2} - 10^{-3} \ll 1.$$

---

[5] При этом остальные ограничения имеют такой же вид, а функционал $F(d) = -\sum_{i,j=1}^{n} d_{ij} \ln\left(d_{ij}/e\right)$. Тогда, согласно экономической интерпретации двойственных множителей Л.В. Канторовича [22]:

$$\beta(C) = \partial F\left(d(C)\right)/\partial C.$$

Из такой интерпретации иногда делают вывод о том, что $\beta$ можно понимать, как цену единицы времени в пути. Чем больше $C$, тем меньше $\beta$.

[6] Отметим, что также как и в [17] из принципа ле Шателье–Самуэльсона [22] следует, что с ростом $\beta$ среднее время в пути $\sum_{i,j=1}^{n} d_{ij}(\beta)T_{ij}\left(d(\beta)\right)$ будет убывать. В связи с этим обстоятельством, а также исходя из соображений размерности, вполне естественно понимать под $\beta$ величину, обратную к характерному (среднему) времени в пути [17] – физическая интерпретация. Собственно, такая интерпретация параметра $\beta$, как правило, и используется в многостадийных моделях (см., например, [1], а также п. 5).

[7] Точнее говоря, из теоремы Т. Куртца [35] следует (с некоторыми, довольно общими, оговорками относительно зависимостей $T_{ij}(d/N)$): если предположить, что делается такой предельный переход, что

$$\exists\ c_{ij}(0) = \lim_{N\to\infty} d_{ij}(0)/N \text{, то } \forall\ t \ge 0\ \exists\ c_{ij}(t) \overset{п.н.}{=} \lim_{N\to\infty} d_{ij}(t)/N$$

(в общем случае пределы существуют не равномерно по времени, но в нашем случае равномерно), где функции $c_{ij}(t)$ (не случайные) удовлетворяют выписанной системе ОДУ.

Для обоснования возможности применения здесь такого рода результата как теорема Тихонова требуется много усилий: во-первых, изначально введенные динамики стохастические (см. также пп. 4, 7, 8 для $T_{ij}$), поэтому, в конечном итоге, нужно все обосновывать именно для них, во-вторых, нам интересна асимптотика по времени первой системы – в медленном времени (то есть нельзя ограничиться классическим случаем: ограниченного отрезка времени), в-третьих, сложный характер оператора, стоящего в правой части второй системы – в быстром времени, не позволяет явно его выписать. Тем не менее, процедура построения этого оператора, которая, в свою, очередь, предполагает переход к пределу (по $\mu \to 0+$ и по числу пользователей транспортной сети (аналог $N(\to\infty)$) – см. п. 7), может быть при необходимости получена из того, что будет далее приведено в пп. 4, 7, 8. Мы не будем здесь подробно описывать, как можно бороться с указанными проблемами, заметим лишь, что в виду “хороших” свойств зависимостей $T_{ij}(c)$ (см. пп. 7, 8), получающихся из приравнивания нулю левой части системы в быстром времени, функция[8]

$$H(c) = \sum_{i,j=1}^{n} c_{ij} \ln(c_{ij}) + \beta\Phi(c), \text{ (следует сравнить с (1))}$$

где $\partial\Phi(c)/\partial c_{ij} = T_{ij}(c)$ (то, что такая функция $\Phi$ существует, показано ниже), будет функцией Ляпунова, и, одновременно, функцией Санова (действием), то есть функцией, которая с точностью до знака и аддитивной постоянной характеризует экспоненциальную концентрацию стационарной меры введенной марковской динамики. Это довольно общий факт, справедливый для широкого класса систем [6, 20]. Поиск равновесной конфигурации $\{c_{ij}^*\}_{i=1,\, j=1}^{n,\, n}$ представляет собой решение задачи

$$\min\left\{ H(c) : c \ge 0, \sum_{j=1}^{n} c_{ij} = l_i, \sum_{i=1}^{n} c_{ij} = \mathrm{w}_j \right\},$$

где $l = L/N$, $\mathrm{w} = W/N$, то есть (1). По $\{c_{ij}^*\}_{i=1,\, j=1}^{n,\, n}$ ($\{d_{ij}^*\}_{i=1,\, j=1}^{n,\, n} = N\{c_{ij}^*\}_{i=1,\, j=1}^{n,\, n}$) уже можно будет определить и равновесные $T_{ij} = T_{ij}(c^*)$.

В дальнейшем, нам будет удобно привести задачу (1) к следующему виду (при помощи метода множителей Лагранжа [14], теоремы фон Неймана о минимаксе [15] и перенормировки $d := d/N$):

$$\min_{\lambda^L, \lambda^W} \max_{\sum_{i,j=1}^{n} d_{ij}=1,\, d_{ij}\ge 0} \left[ -\sum_{i,j=1}^{n} d_{ij} \ln d_{ij} - \beta \sum_{i=1,\, j=1}^{n,\, n} d_{ij} T_{ij} + \sum_{i=1}^{n} \lambda_i^L \left( l_i - \sum_{j=1}^{n} d_{ij} \right) + \sum_{i=1}^{n} \lambda_j^W \left( \mathrm{w}_j - \sum_{i=1}^{n} d_{ij} \right) \right], \quad (2)$$

[8] Полученная функция будет выпуклой. Далее это будет проясняться.

где, напомним, $l = L/N$, $\mathrm{w} = W/N$.[9] Вместо того, чтобы подставлять сюда зависимость $T_{ij}(d)$, мы используем следующий трюк. В пп. 6–8 задача поиска зависимости $T_{ij}(d)$ будет сведена к задаче вида: $\min_{t \ge \bar{t}} \left\{ \beta \langle \bar{f}, t - \bar{t} \rangle - \beta \sum_{i=1,\, j=1}^{n,\, n} d_{ij} T_{ij}(t) \right\} \stackrel{def}{=} -\Phi(d)$, где $\bar{t}$ и $\bar{f}$ – известные векторы (входные параметры модели) временных затрат в пути на ребрах графа транспортной сети и максимальных пропускных способностей ребер графа транспортной сети, $T_{ij}(t)$ – длина кратчайшего пути из района $i$ в район $j$ (на графе транспортной сети имеется много вершин – намного больше числа районов, но мы считаем, что в каждом районе есть только одна вершина транспортного графа, являющаяся источником/стоком для пользователей сети, именно между такими представителями вершин районов $i$ и $j$ берется кратчайшие расстояние), если веса ребер графа транспортной сети задаются вектором $t$. Не стоит путать $T_{ij}(t)$ с искомой зависимостью $T_{ij}(d)$ – это разные функциональные зависимости. Решив указанную задачу $t^*(d)$ минимизации, считая $d$ параметрами, мы получили бы искомую зависимость $T_{ij}(d) := T_{ij}(t^*(d))$. Но явно это нельзя сделать в типичных ситуациях, да и не нужно, потому что в итоге все равно необходимо работать с потенциалом $\Phi(d)$. Поэтому предлагается ввести в задачу (1) подзадачу и добавить слагаемое $\min_{\lambda^L, \lambda^W} \max_{\sum_{i,j=1}^{n} d_{ij} = 1,\, d_{ij} \ge 0} \min_{t \ge \bar{t}} \left[ \cdot + \beta \langle \bar{f}, t - \bar{t} \rangle \right]$. Решение такой задачи сразу даст все, что нужно. Это следует из формулы Демьянова–Данскина [14]. Пункт 9 посвящен упрощению только что сформулированной конструкции. Более полное обоснование и дальнейшее развитие можно отследить по работам [36–38].

[9] Обратим внимание, что мы берем максимум в (2) при дополнительном ограничении $\sum_{i,j=1}^{n} d_{ij} = 1$, которого изначально не было. Однако легко показать, что это ограничения является следствием системы ограничений (А), точнее следствием одновременно двух подсистем в (А) – отвечающих $l$ и $w$. Как следствие, можно считать, что $\sum_{i=1}^{n} \lambda_i^L = 0$, $\sum_{j=1}^{n} \lambda_j^W = 0$. Отметим также, что задача (2) может быть упрощена, поскольку внутренний максимум явно находится, однако мы отложим соответствующие выкладки до п. 9. Здесь же мы рассмотрим случай, когда мы заносим в функционал с помощью множителей Лагранжа только часть ограничений (соответствующих $l$ или $w$). В таком случае, вычисляя внутренний максимум по $d$, получим, соответственно, задачи:

$$\min_{\lambda^W} \left[ \sum_{i=1}^{n} l_i \ln \left[ \sum_{j=1}^{n} \exp\left(-\lambda_j^W - \beta T_{ij}\right) \right] + \sum_{j=1}^{n} \lambda_j^W \mathrm{w}_j \right]; \quad \min_{\lambda^L} \left[ \sum_{j=1}^{n} w_j \ln \left[ \sum_{i=1}^{n} \exp\left(-\lambda_i^L - \beta T_{ij}\right) \right] + \sum_{i=1}^{n} \lambda_i^L l_i \right],$$

где, соответственно

$$d_{ij} = l_i \exp\left(-\lambda_j^W - \beta T_{ij}\right) \left[ \sum_{k=1}^{n} \exp\left(-\lambda_k^W - \beta T_{ik}\right) \right]^{-1}; \quad d_{ij} = w_j \exp\left(-\lambda_i^L - \beta T_{ij}\right) \left[ \sum_{k=1}^{n} \exp\left(-\lambda_k^L - \beta T_{kj}\right) \right]^{-1}.$$

Отсюда можно усмотреть интерпретацию двойственных множителей как соответствующих “потенциалов притяжения/отталкивания районов” [17]. К этому мы еще вернемся в замечании п. 4. В принципе далее удобнее было бы работать с одной из этих задач, а не с (2), однако для сохранения симметрии, мы оставляем задачу в форме (2).

Подобно тому, как мы рассматривали в этом пункте трудовые корреспонденции (в утренние и вечерние часы более 70% корреспонденций по Москве и области именно такие), можно рассматривать перемещения, например, из дома к местам учебы, отдыха, в магазины и т.п. (по хорошему, еще надо было учитывать перемещения типа работа–магазин–детский_сад–дом) – рассмотрение всего этого вкупе приведет также к задаче (1). Только будет больше типов корреспонденций $d$ : помимо пары районов, еще нужно будет учитывать тип корреспонденции [глава 2, 17]. Все это следует из того, что инвариантной мерой динамики с несколькими типами корреспонденций по-прежнему будет прямое произведение пуассоновских мер. Важно отметить при этом, что в таком контексте (равновесные) $T_{ij}$ будут определяться парами районов, а не типом передвижения, то есть с точки зрения последующего изложения это означает, что ничего по-сути не поменяется. Другое дело, когда мы рассматриваем разного типа пользователей транспортной сети, например: имеющих личный автомобиль и не имеющих личный автомобиль. Первые могут им воспользоваться равно, как и общественным транспортом, а вторые нет. И на рассматриваемых масштабах времени пользователи могут менять свой тип. То есть время в пути может для разных типов пользователей быть различным [глава 2, 17]. Считая, подобно тому как мы делали раньше, что желание пользователей корреспонденции[10] $(i, j)$ сети сменить свой тип (вероятность в единицу времени) есть

$$\tilde{\lambda} \exp\Big( \underbrace{\tilde{R}\left(T_{ij}^{old}\right)}_{\substack{\text{суммарные затраты} \\ \text{до смены типа}}} - \underbrace{\tilde{R}\left(T_{ij}^{new}\right)}_{\substack{\text{суммарные затраты} \\ \text{после смены типа}}} \Big), \text{ где } \tilde{R}(T) = \beta T ,$$

и учитывая в “обменах” тип пользователя (будет больше типов корреспонденций $d$, но “меняются местами работы” только пользователи одного типа), можно показать, что все это вкупе приведет также к задаче типа (1), и все последующие рассуждения распространяются на этот случай. Но с оговоркой, что в п. 8 расщепление потоков надо делать только для тех пользователей, для которых имеется возможность использовать личный транспорт. При этом часть пользователей (в зависимости от текущего $d$) распределяются только по сети общественного транспорта. Несложно понять, что это ничего принципиально не изменит. Фактически, в этом абзаце мы описали, как сделать из трехстадийной модели полноценную четырехстадийную модель, которыми обычно и пользуются на практике. Чтобы не делать изложение излишне громоздким, далее мы не учитываем нюансы, описанные в этом абзаце.

[10] Здесь, конечно, надо учитывать не только трудовые миграции, но и все остальные, поскольку для заметной части жителей Москвы и области решение о покупке автомобиля напрямую связано с желанием ездить на нем в основном только на дачу. В этом месте возникает необходимость моделирования (учета) перемещений пользователей сети не только в рамках установленного диапазона времени в течение типичного дня, но и в целом все возможные перемещения. Это обстоятельство вынуждает использовать здесь различного рода эвристики, дабы не отказываться от важного предположения: “в рамках установленного диапазона времени в течение типичного дня”.

# 4. Модель равновесного распределения потоков Бэкмана

Следуя книге [6], опишем наиболее популярную на протяжении более чем полувека модель равновесного распределения потоков Бэкмана [2–6].

Пусть транспортная сеть города представлена ориентированным графом $\Gamma = (V, E)$, где $V$ – узлы сети (вершины), $E \subset V \times V$ – дуги сети (рёбра графа). В современных моделях равновесного распределения потоков в крупном мегаполисе число узлов графа транспортной сети обычно выбирают порядка $|V| \sim 10^4 - 10^5$. Число ребер $|E|$ получается в три четыре раза больше. Пусть $W \subseteq \{w = (i, j) : i, j \in V\}$ – множество корреспонденций, т.е. возможных пар «исходный пункт» – «цель поездки» ($|W|$ по порядку величины может быть $n^2$); $p = \{v_1, v_2, ..., v_m\}$ – путь из $v_1$ в $v_m$, если $(v_k, v_{k+1}) \in E$, $k = 1, ..., m-1$, $m > 1$; $P_w$ – множество путей, отвечающих корреспонденции $w \in W$, то есть если $w = (i, j)$, то $P_w$ – множество путей, начинающихся в вершине $i$ и заканчивающихся в $j$; $P = \bigcup_{w \in W} P_w$ – совокупность всех путей в сети $\Gamma$ (число “разумных” маршрутов $|P|$, которые потенциально могут использоваться, обычно растет с ростом числа узлов сети не быстрее чем $O\left(|V|^3\right)$, однако теоретически может быть экспоненциально большим); $x_p$ [автомобилей/час] – величина потока по пути $p$, $x = \{x_p : p \in P\}$; $f_e$ [автомобилей/час] – величина потока по дуге $e$:

$$f_e(x) = \sum_{p \in P} \delta_{ep} x_p \text{, где } \delta_{ep} = \begin{cases} 1, & e \in p \\ 0, & e \notin p \end{cases};$$

$\tau_e(f_e)$ – удельные затраты на проезд по дуге $e$. Как правило, предполагают, что это – (строго) возрастающие, гладкие функции от $f_e$ (в конце этого раздела нам потребуется еще и выпуклость). Точнее говоря, под $\tau_e(f_e)$ правильнее понимать представление пользователей транспортной сети об оценке собственных затрат (обычно временных в случае личного транспорта и комфортности пути (с учетом времени в пути) в случае общественного транспорта) при прохождении дуги $e$, если поток желающих оказаться на этой дуге будет $f_e$.

Зависимость $\tau_e(f_e)$ можно попробовать и вывести, например, из следующих соображений [6] (вариация на тему модели Бобкова–Буслаева–Танака). Рассматривается одна полоса длины $L$ и транспортный поток, характеризующийся максимальной скоростью $v_{\max}$ и следующей зависимостью безопасного расстояния (“комфортного” расстояния до впереди идущего транспортного средства) от скорости: $d(v) = l + \tilde{\tau} v + c v^2$, где $l$ – средняя длина автомобиля в “стоячей” пробке (эта длина немного больше средней “физической” длины автомобиля $\approx 6.5$ м), $\tilde{\tau}$ – время реакции водителей (эксперименты показывают, что для европейских водителей эта величина обычно равна одной секунде, для российских водителей она, как правило, не превышает полсекунды), $c$ – характеризует коэффициент трения шин о поверхность дороги, поскольку слагаемое $cv^2$ характеризует в $d(v)$ тормозной путь. Действительно, пока водитель среагирует на ситуацию он в среднем проедет (не изменяя своей скорости) путь $\tilde{\tau}v$. Когда уже реакция произошла, водитель начинает тормозить и кинетическая энергия автомобиля $mv^2/2$ должна быть “погашена” работой силой трения на участке $h$ тормозного пути $\mu mgh$ ($\mu$ – коэффициент трения, $g$ – ускорение свободного падения). Отсюда можно найти $c = 1/(2\mu g)$. Имея зависимость $d(v)$, можно ввести зависимость $\rho(v) = 1/d(v)$, которая порождает зависимость потока от скорости $f(v) = v\rho(v)$. Используя то, что $\tau(f(v)) = L/v$ и $0 \le v \le v_{\max}$ можно явно выписать искомую зависимость $\tau(f)$.

Рассмотрим теперь $G_p(x)$ – затраты временные или финансовые на проезд по пути $p$. Естественно считать, что $G_p(x) = \sum_{e\in E} \tau_e(f_e(x))\delta_{ep}$. В приложениях часто требуется учитывать также затраты на прохождения вершин графа, которые могут зависеть от величин всех потоков через рассматриваемую вершину.

Пусть также известно, сколько перемещений в единицу времени $d_w$ осуществляется согласно корреспонденции $w \in W$. Тогда вектор $x$, характеризующий распределение потоков, должен лежать в допустимом множестве:

$$X = \left\{ x \ge 0 : \sum_{p\in P_w} x_p = d_w,\ w \in W \right\}.$$

Это множество может иметь и более сложный вид, если дополнительно учитывать, например, конечность пропускных способностей рёбер (ограничения сверху на $f_e$).

Рассмотрим игру, в которой каждому элементу $w \in W$ соответствует свой, достаточно большой ($d_w \gg 1$), набор однотипных “игроков”, осуществляющих передвижение согласно корреспонденции $w$. Чистыми стратегиями игрока служат пути, а выигрышем – величина $-G_p(x)$. Игрок “выбирает” путь следования $p \in P_w$, при этом, делая выбор, он пренебрегает тем, что от его выбора также “немного” зависят $|P_w|$ компонент вектора $x$ и, следовательно, сам выигрыш $-G_p(x)$. Можно показать, что отыскание равновесия Нэша–Вардропа $x^* \in X$ (макро описание равновесия) равносильно решению задачи нелинейной комплементарности (принцип Вардропа):

$$\textit{для любых } w \in W,\ p \in P_w \textit{ выполняется } x_p^* \cdot \left( G_p(x^*) - \min_{q\in P_w} G_q(x^*) \right) = 0.$$

Действительно допустим, что реализовалось какое-то другое равновесие $\tilde{x}^* \in X$, которое не удовлетворяет этому условию. Покажем, что тогда найдется водитель, которому выгодно поменять свой маршрут следования. Действительно, тогда существуют такие $\tilde{w} \in W,\ \tilde{p} \in P_{\tilde{w}}$, что $\tilde{x}_{\tilde{p}}^* \cdot \left( G_{\tilde{p}}\left(\tilde{x}^*\right) - \min_{q \in P_{\tilde{w}}} G_q\left(\tilde{x}^*\right) \right) > 0$.

Каждый водитель (множество таких водителей не пусто $\tilde{x}_{\tilde{p}}^* > 0$), принадлежащий корреспонденции $\tilde{w} \in W$, и использующий путь $\tilde{p} \in P_{\tilde{w}}$, действует не разумно, поскольку существует такой путь такой путь $\tilde{q} \in P_{\tilde{w}},\ \tilde{q} \neq \tilde{p}$, что $G_{\tilde{q}}\left(\tilde{x}^*\right) = \min_{q \in P_{\tilde{w}}} G_q\left(\tilde{x}^*\right)$. Этот путь $\tilde{q}$ более выгоден, чем $\tilde{p}$. Аналогично показывается, что при $x^* \in X$ никому из водителей уже не выгодно отклоняться от своих стратегий. Но это по определению и называется равновесием Нэша, которое ввел в своей диссертации в конце 40-х годов XX века Джон Нэш, получивший именно за эту концепцию в 1994 г. нобелевскую премию по экономике [39]. Мы также добавляем фамилию Дж.Г. Вардропа, которой чуть позже Дж. Нэша привнес к этой концепции условие “конкурентного рынка”: игрок, принимающий решение, пренебрегает тем, что его решение сколько-нибудь значительно поменяет ситуацию на “рынке”. Когда игроков двое, трое (ситуации, рассматриваемые Дж. Нэшем), то, очевидно, что так делать нельзя. Но когда игроков (водителей) десятки и сотни тысяч… Вся эта конструкция неявно предполагает, что $x_p^* > 0 \Rightarrow x_p^* \gg 1$. Поэтому, не боясь сильно ошибиться, можно искать решение задачи нелинейной комплементарности, не предполагая целочисленности компонент вектора $x^* \in X$. Такая релаксация изначально целочисленной задачи заметно упрощает её с вычислительной точки зрения!

Хотя мы и смогли выписать условие равновесия в виде задачи нелинейной комплементарности, это не сильно продвинуло нас в понимании того, как его находить. Пытаться честно решить задачу в таком виде – вычислительно бесперспективная задача. С другой стороны современные вычислительные методы позволяют эффективно решать задачи выпуклой оптимизации. Постараемся свести нашу задачу к таковой.

Для этого, прежде всего, заметим, что рассматриваемая нами игра принадлежит к классу, так называемых, потенциальных игр. В нашем случае это означает, что существует такая функция

$$\Psi(x) = \sum_{e \in E} \int_0^{\sum_{p \in P} \delta_{ep} x_p} \tau_e(z)\,dz = \sum_{e \in E} \sigma_e\left(f_e(x)\right),$$

где $\sigma_e\left(f_e\right) = \int_0^{f_e(x)} \tau_e(z)\,dz$, что $\partial\Psi(x)/\partial x_p = G_p(x)$ для любого $p \in P$. Таким образом, мы имеем дело с потенциальной игрой. Оказывается, что $x^* \in X$ – равновесие Нэша–Вардропа тогда и только тогда, когда оно доставляет минимум $\Psi(x)$ на множестве $X$.

Действительно, предположим, что $x^* \in X$ – точка минимума. Тогда, в частности, для любых $w \in W$, $p, q \in P_w$ ($x_p^* > 0$) и достаточно маленького $\delta x_p > 0$ выполняется:

$$-\frac{\partial \Psi(x^*)}{\partial x_p}\delta x_p + \frac{\partial \Psi(x^*)}{\partial x_q}\delta x_p \ge 0 .$$

Иначе, заменив $x^*$ на

$$\breve{x}^* = x^* + \left( \underbrace{\overbrace{0,...,0,-\delta x_p}^{p},0,...,0,\delta x_p}_{q},0,...,0 \right) \in X ,$$

мы пришли бы к вектору $\breve{x}^*$, доставляющему меньшее значение $\Psi(x)$ на множестве $X$ :

$$\Psi(\breve{x}^*) \approx \Psi(x^*) - \frac{\partial \Psi(x^*)}{\partial x_p}\delta x_p + \frac{\partial \Psi(x^*)}{\partial x_q}\delta x_p < \Psi(x^*).$$

Вспоминая, что $\partial \Psi(x)/\partial x_p = G_p(x)$, и учитывая, что $q$ можно выбирать произвольно из множества $P_w$, получаем:

*для любых* $w \in W$, $p \in P_w$, *если* $x_p^* > 0$, *то выполняется* $\min\limits_{q \in P_w} G_q(x^*) \ge G_p(x^*)$.

Но это и есть по-другому записанное условие нелинейной комплементарности. Строго говоря, мы показали сейчас только то, что точка минимума $\Psi(x)$ на множестве $X$ будет равновесием Нэша–Вардропа. Аналогично рассуждая, можно показать и обратное: равновесие Нэша–Вардропа доставляет минимум $\Psi(x)$ на множестве $X$ . Этот минимум можно искать, например, с помощью метода условного градиента [18, 40] (Франк–Вульф).

**Теорема 1 [3–6].** *Вектор* $x^*$ *будет равновесием Нэша–Вардропа тогда и только тогда, когда*

$$x \in \operatorname{Arg}\min_x \left[ \Psi(f(x)) = \sum_{e \in E} \sigma_e(f_e(x)) :\ f = \Theta x,\ x \in X \right]. \quad (3)$$

*Если преобразование* $G(\cdot)$ *строго монотонное, то равновесие* $x$ *единственно. Если* $\tau_e'(\cdot) > 0$, *то равновесный вектор распределения потоков по ребрам* $f$ *– единственный* (*это еще не гарантирует единственность вектора распределения потоков по путям* $x$ *[6]*).

Проинтерпретируем, следуя [6, 9, 21, 41], эволюционным образом полученное равновесие, попутно отвечая на вопрос (поскольку задача (4) ниже имеет единственное решение, детали см. в [21, 42]): какому из равновесий стоит отдать предпочтение, в случае неединственности? Опишем марковскую логит динамику (также говорят гиббсовскую динамику) в повторяющейся игре загрузки графа транспортной сети [9]. Пусть каждой корреспонденции отвечает $d_w M$ агентов ($M \gg 1$), $x := x/M$, $f := f/M$, $\tau_e(f_e) := \tau_e(f_e/M)$.

Пусть имеется $tN$ шагов ($N \gg 1$). Пусть $k$-й агент, принадлежащий корреспонденции $w \in W$, независимо от остальных на шаге $m+1$ с вероятностью с вероятность $1-\lambda/N$ выбирает путь $p^{k,m}$, который использовал на шаге $m = 0,...,tN$, а с вероятностью $\lambda/N$ ($\lambda > 0$) решает "поменять" путь, и выбирает (возможно тот же самый) зашумленный кратчайший путь

$$p^{k,m+1} = \arg\max_{q \in P_w} \left\{ -G_q\left(x^m\right) + \xi_q^{k,m+1} \right\},$$

где независимые случайные величины $\xi_q^{k,m+1}$, имеют одинаковое двойное экспоненциальное распределение, также называемое распределением Гумбеля[11] [9, 41]:

$$P\left(\xi_q^{k,m+1} < \zeta\right) = \exp\left\{-e^{-\zeta/T-E}\right\},\ T > 0.$$

Отметим, что

$$P\left(p^{k,m+1} = p \,\middle|\, \text{агент решил "поменять" путь}\right) = \frac{\exp\left(-G_p\left(x^m\right)/T\right)}{\sum_{q \in P_w} \exp\left(-G_q\left(x^m\right)/T\right)}.$$

Кроме того, если взять $E \approx 0.5772$ – константа Эйлера, то

$$M\left[\xi_q^{k,m+1}\right] = 0,\ D\left[\xi_q^{k,m+1}\right] = T^2\pi^2/6.$$

Такая динамика отражает ограниченную рациональность агентов (см. п. 5), и часто используется в популяционной теории игр [9] и теории дискретного выбора [41]. Оказывается эта марковская динамика в пределе $N \to \infty$ превращается в марковскую логит динамику в непрерывном времени (вырождающуюся при $T \to 0+$ в динамику наилучших ответов [9] – последующие рассуждения, в частности, формулы (4), (5), допускают переход к пределу $T \to 0+$). Марковская логит динамика в непрерывном времени допускает два предельных перехода (обоснование перестановочности этих пределов см. в [34, 35]): $t \to \infty$, $M \to \infty$ или $M \to \infty$, $t \to \infty$. При первом порядке переходов мы сначала ($t \to \infty$) согласно эргодической теореме для марковских процессов (в нашем случае марковский процесс – модель стохастической химической кинетики с унарными реакциями в условиях детального баланса [21, 34]) приходим к финальной (=стационарной) вероятностной мере, имеющей в основе мультиномиальное распределение. С ростом числа агентов ($M \to \infty$) эта мера

$$\sim \exp\left(-\frac{M}{T}\cdot\left(\Psi_T\left(x\right) + o(1)\right)\right)$$

---

[11] Распределение Гумбеля можно объяснить исходя из идемпотентного аналога центральной предельной теоремы (вместо суммы случайных величин – максимум) для независимых случайных величин с экспоненциальным и более быстро убывающим правым хвостом [43]. Распределение Гумбеля возникает в данном контексте, например, если при принятии решения водитель собирает информацию с большого числа разных (независимых) зашумленных источников, ориентируясь на худшие прогнозы по каждому из путей.

экспоненциально концентрируется около наиболее вероятного состояния, поиск которого сводится к решению энтропийно регуляризованной задачи (3) (как численно решать эту задачу описано в работе [44], см. также [32])

$$\Psi_T(x) = \Psi(f(x)) + T\sum_{w\in W}\sum_{p\in P_w} x_p \ln x_p \to \min_{\substack{f(x)=\Theta x\\ x\in X}}. \tag{4}$$

Функционал в этой задаче оптимизации с точностью до потенцирования и мультипликативных и аддитивных констант соответствует исследуемой стационарной мере – то есть это функционал Санова [21, 34]. При обратном порядке предельных переходов, мы сначала ($M\to\infty$) осуществляем, так называемый, канонический скейлинг [34, 35], приводящий к детерминированной кинетической динамике, описываемой СОДУ на $x$

$$\frac{dx_p}{dt} = d_w \frac{\exp(-G_p(x)/T)}{\sum_{l\in P_w}\exp(-G_l(x)/T)} - x_p,\ p\in P_w,\ w\in W, \tag{5}$$

а затем ($t\to\infty$) ищем аттрактор получившейся СОДУ. Глобальным аттрактором оказывается неподвижная точка, которая определяется решением задачи (4). Более того, функционал $\Psi_T(x)$ является функцией Ляпунова полученной кинетической динамики (5) (то есть является функционалом Больцмана). Последнее утверждение – достаточно общий факт (функционал Санова, является функционалом Больцмана), верный при намного более общих условиях (см. [34] и цитированную там литературу).

Хотелось бы подчеркнуть, что рассматриваемая выше “игра” – потенциальная (это общий факт для игр загрузок [9]; Розенталь 1973, Мондерер–Шэпли 1996), поэтому из общих результатов эволюционной теории игр [9, 45], следует, что любые разумные содержательно интерпретируемы (суб-)градиентные спуски приводят к равновесию Нэша–Вардропа (или его стохастическому варианту, даваемому решением задачи (4)). В частности, в [6, 46, 47] содержательно интерпретируется быстро сходящаяся динамика, связанная с методом зеркального спуска, которая порождается имитационной логит динамикой [9]. Хотелось бы также обратить внимание на эволюционную интерпретацию парадокса Браесса: когда не эффективное по Парето, единственное равновесие Нэша–Вардропа в специально сконструированной транспортной сети является, тем не менее, эволюционно устойчивым [6, 9].

Нетривиальным является следующее наблюдение. Если рассмотреть энтропийно-регуляризованный функционал $\Psi_T(x)$ и взять предел (см. пп. 6, 7)

$$\tau_e^{\mu}(f_e) \xrightarrow[\mu\to 0+]{} \begin{cases} \bar{t}_e, & 0\le f_e < \bar{f}_e \\ [\bar{t}_e,\infty), & f_e = \bar{f}_e \end{cases},$$

$$d\tau_e^{\mu}(f_e)/df_e \xrightarrow[\mu\to 0+]{} 0,\ 0\le f_e < \bar{f}_e,$$

то переход к двойственной задаче (для задачи минимизации этого функционала на множестве $X$ ) дает стохастический вариант модели стабильной динамики (см. п. 6), который используется в стохастическом варианте трехстадийной модели стабильной динамики (см. п. 10).

**Замечание 1 (“облачная модель” расчета матрицы корреспонденций [21]).** В контексте написанного выше полезно отметить другой способ обоснования энтропийной модели расчета матрицы корреспонденций из п. 3. Предположим, что все вершины, отвечающие источникам корреспонденций, соединены ребрами с одной вспомогательной вершиной (облако № 1). Аналогично, все вершины, отвечающие стокам корреспонденций, соединены ребрами с другой вспомогательной вершиной (облако № 2). Припишем всем новым ребрам постоянные веса. И проинтерпретируем веса ребер, отвечающих источникам $\lambda_i^L$, например, как средние затраты на проживание (в единицу времени, скажем, в день) в этом источнике (районе), а веса ребер, отвечающих стокам $\lambda_j^W$, как уровень средней заработной платы со знаком минус (в единицу времени) в этом стоке (районе), если изучаем трудовые корреспонденции. Будем следить за системой в медленном времени, то есть будем считать, что равновесное распределение потоков по путям стационарно. Поскольку речь идет о равновесном распределении потоков, то нет необходимости говорить о затратах на путях или ребрах детализированного транспортного графа, достаточно говорить только затратах (в единицу времени), отвечающих той или иной корреспонденции. Таким образом, у нас есть взвешенный транспортный граф с одним источником (облако 1) и одним стоком (облако 2). Все вершины этого графа, кроме двух вспомогательных (облаков), соответствуют районам в модели расчета матрицы корреспонденций из п. 3. Все ребра этого графа имеют стационарные (не меняющиеся и не зависящие от текущих корреспонденций) веса $\left\{T_{ij};\lambda_i^L;\lambda_j^W\right\}$. Если рассмотреть естественную в данном контексте логит динамику с $T = 1/\beta$ (здесь полезно напомнить, что согласно п. 3 $\beta$ обратно пропорционально средним затратам, а $T$ имеет как раз физическую размерность затрат), описанную выше, то поиск равновесия рассматриваемой макросистемы приводит (в прошкалированных переменных) к задаче, сильно похожей на задачу (2) из п. 3

$$\max_{\sum_{i,j=1}^{n} d_{ij}=1,\ d_{ij}\ge 0}\left[-\sum_{i,j=1}^{n} d_{ij}\ln d_{ij} - \beta\sum_{i=1,\ j=1}^{n,n} d_{ij}T_{ij} - \beta\sum_{i=1}^{n}\left(\lambda_i^L\sum_{j=1}^{n} d_{ij}\right) - \beta\sum_{i=1}^{n}\left(\lambda_j^W\sum_{i=1}^{n} d_{ij}\right)\right].$$

Разница состоит в том, что здесь мы не оптимизируем по $2n$ двойственным множителям $\lambda^L$, $\lambda^W$. Более того, мы их и не интерпретируем здесь как двойственные множители, поскольку мы их ввели на этапе взвешивания ребер графа. Тем не менее, значения этих переменных, как правило, не откуда брать. Тем более, что приведенная выше (наивная) интерпретация вряд ли может всерьез рассматриваться, как способ определения этих пара-метров исходя из данных статистики. Более правильно понимать $\lambda_i^L$, $\lambda_j^W$ как потенциалы притяжения/отталкивания районов, включающиеся в себя плату за жилье и зарплату, но включающие также и многое другое, что сложно описать количественно. И здесь как раз помогает информация об источниках и стоках, содержащаяся в $2n$ уравнениях из формулы (**А**) п. 3. Таким образом, мы приходим ровно к той же самой задаче (2) с той лишь разницей, что мы получили дополнительную интерпретацию двойственных множителей в задаче (2). При этом двойственные множители в задаче (2) равны (с точностью до мультипликативного фактора $\beta$) введенным здесь потенциалам притяжения районов. Несложно распространить на изложенную здесь модель написанное в п. 3 по поводу того, как с помощью разделения времен можно учитывать обратную связь: перераспределение потоков по путям изменяется (в быстром времени) при изменении корреспонденций, а также распространить на то, что написано в самом конце п. 3. Нам представляется такой способ рассуждения даже более привлекательным, чем способ, описанный в п. 3, и основанный на "обменах". И связано это с тем, что для получения равновесия в многостадийной модели, мы можем рассмотреть всего одну (общую) логит динамику, в которой с малой интенсивностью (в медленном времени) происходят переходы, описанные в этом замечании (жители города меняют места жительства, работы), а с высокой интенсивностью (в быстром времени, изо дня в день) жители города перераспределяются по путям (в зависимости от текущих корреспонденций, подстраиваясь под корреспонденции) – это как раз и было описано непосредственно перед замечанием. Другая причина – большая вариативность модели построенной в этом замечании. Нам представляется очень плодотворной и перспективной идея перенесения имеющейся информации об исследуемой системе из обременительных законов сохранения динамики, описывающей эволюцию этой системы, в саму динамику путем введения дополнительных естественно интерпретируемых параметров. При таком подходе становится возможным, например, учитывать в моделях и рост транспортной сети. Другими словами, при таком подходе, например, можно естественным образом рассматривать также и ситуацию, когда число пользователей транспортной сетью меняется со временем (медленно).

В заключение отметим, что если штрафовать (назначать платы) за использование различных стратегий (маршрутов) по правилу

$$\bar{G}_p(x) = G_p(x) + \underbrace{\sum_{w\in W}\sum_{q\in P_w} x_q \frac{\partial G_q(x)}{\partial x_p}}_{\textit{штраф}},\ p \in P,$$

то найдется такая функция $\bar{\Psi}(x) = \sum_{w\in W}\sum_{p\in P_w} x_p G_p(x)$, что $\partial\bar{\Psi}(x)/\partial x_p = \bar{G}_p(x)$. Поэтому из сказанного выше в этом пункте будет следовать, что возникающее в такой управляемой транспортной сети равновесие Нэша–Вардропа будет (единственным) глобально устойчи-вым и соответствовать социальному оптимуму в изначальной транспортной сети [8]. Чтобы в этом убедиться, достаточно (в виду линейной связи $f = \Theta x$) проверить строгую выпуклость функции $\bar{\Psi}(x) = \sum_{w \in W} \sum_{p \in P_w} x_p G_p(x) = \sum_{e \in E} f_e(x) \tau_e(f_e(x))$, для чего достаточно выпуклости функций $\tau_e(f_e)$, $e \in E$. Все это хорошо соответствует механизму Викри–Кларка–Гроуса [48] (VCG mechanism) – штраф (плата) за использование маршрута новым пользователем равен дополнительным потерям, которые понесут из-за этого все остальные пользователи. Однако важно сделать две оговорки. Во-первых, все это хотя и можно попытаться практически осуществить (например, собирая транспортные налоги исходя из трековой информации, которую в перспективе можно будет иметь о каждом автомобиле), но довольно сложным оказывается механизм. Платы взимаются не за проезд по ребру, как хотелось бы, а именно за выбор (проезд) маршрута. Кроме того, плата является функцией состояния транспортной системы $x$, которое, в отличие от $f$, не наблюдаемо. Во-вторых, взимая платы за проезд, мы с одной стороны приводим систему в социальный оптимум, а с другой стороны для достижения этой цели вынуждены собирать с участников движения налоги. К сожалению, их размер может оказаться внушительным, и это уже необходимо учитывать с точки зрения расщепления участников движения по выбору типа передвижения. Относительно второй проблемы – готовых решений нам не известно. Это известная проблема в современном разделе теории игр: mechanism design [48]. А вот по первой проблеме (адаптивное) решение есть [8]:

$$\bar{\tau}_e(f_e) = \tau_e(f_e) + \underbrace{f_e \tau_e'(f_e)}_{штраф}, \; e \in E.$$

Легко проверить, что это приведет к указанному выше пересчету $G_p(x) \to \bar{G}_p(x)$. Далее, если мы зафиксируем (знаем) социальный оптимум $f^{opt}$, то платы за проезд можно выбирать постоянными $f_e^{opt} \tau_e'(f_e^{opt})$, $e \in E$. Все сказанное выше об устойчивости останется в силе (без всяких дополнительных предположений о выпуклости функций $\tau_e(f_e)$, $e \in E$) с одной лишь оговоркой, что транспортная система должна поддерживаться при заданных корреспонденциях $d_w$ и функциях затрат $\tau_e(f_e)$. В противном случае, возникающее равновесие уже может не соответствовать социальному оптимуму.

# 5. Краткий обзор подходов к построению многостадийных моделей транспортных потоков, с моделью типа Бэкмана в качестве модели равновесного распределения потоков

Так называемая «4-стадийная» модель является наиболее употребительной методологией моделирования транспортных систем городов и агломераций (см., например, [1]). В рамках данной модели производится поиск равновесия спроса и предложения на поездки. При этом рассматриваются в единой совокупности модели генерации трафика, его распределения по типам передвижения и дальнейшее распределение по маршрутам.

Методология, как уже описывалось ранее, включает последовательное выполнение четырех этапов, последние три из которых закольцовываются для получения самосогласованных результатов. Исходными данными для модели являются: граф дорожной сети и сети общественного транспорта с заданными функциями издержек и других определяющих параметров, разделение города на транспортные зоны и параметры этих зон (например, количество рабочих мест или мест жительства). На выходе модели выдаются: оценка матрицы корреспонденций для каждого типа передвижения, загрузка элементов сети (например, конкретной дороги) и издержки, соответствующие данному уровню загрузки. Понятие «тип передвижения», используемое выше, является обобщением понятия тип транспорта и обозначает последовательность (или просто множество) используемых типов транспорта. Например, типом передвижения может считаться поездка на общественном транспорте (не важно, автобусе или троллейбусе или на них обоих поочередно с пересадками) или же использование схемы park-and-ride. Уровень детализации при этом определяется самим модельером. Обычно выбирается или модель с тремя типами передвижений (общественный и личный транспорт, пешие прогулки) или с детализацией до типа транспорта (автомобиль, наземный общественный транспорт, метро, пешие прогулки, и различные комбинации, перечисленные ранее). При этом стараются не учитывать типы передвижений, которые используются очень редко.

Структурно, расчет модели можно описать следующим образом:

1. На первом шаге из исходных данных о транспортных зонах получают вектора отправления и прибытия для транспортных зон, т.е. в наших обозначениях $\{L_i\}_{i=1}^n$ и $\{W_j\}_{j=1}^n$.
2. Рассчитывается первая оценка матрицы корреспонденций (обычно, с помощью гравитационной модели [6, 18]).
3. Рассчитывается расщепление корреспонденций по типам передвижений.
4. Для каждого типа передвижений рассчитывается распределение потоков по маршрутам.

5. Получается оценка матрицы издержек корреспонденций и вектор загрузки сети.
6. Проверяется критерий остановки.
7. Если критерий выполнен, решение получено.
8. Если критерий не выполнен, вернуться на шаг 2 с переоцененной матрицей издержек.

Первый этап модели мы разбирать не будем, так как он является достаточно обособленным от других в том смысле, что он не входит в итерационную часть метода и практически «бесплатен» с точки зрения сложности операций. Отметим лишь, что вектора отправления и прибытия рассчитываются из параметров зон с помощью простейших регрессионных моделей.

Рассмотрим модели, лежащие в основе этого алгоритма более подробно. Данный алгоритм – итерационный, после «прогонки» очередной итерации мы возвращаемся на второй «шаг» схемы. На $m$-й итерации на 2-м шаге формируется $m$-я оценка матрицы корреспонденций. Для её построения необходимо знание о матрице издержек в сети (т.е. знание о стоимости проезда из каждой зоны в каждую). На всех шагах, кроме первого, данные значения получаются из предыдущей ($m$–1)-й итерации. На первой итерации используются издержки, соответствующие незагруженной сети или почти любая другая, разумная, оценка матрицы издержек. Если в сети для пользователей доступно несколько типов передвижений, то в качестве оценки издержек берутся или средние издержки (время в пути) для конкретной корреспонденции или же (если пользователи распределены «равновесно») значения издержек, например, для личного транспорта, так как в «равновесии» издержки зачастую (исключая некоторые особые случаи) должны быть равны для различных, используемых, типов передвижений.

Как правило, оценки матрицы корреспонденций строятся согласно гравитационной или энтропийной модели [6, 18].

Посмотрим на то, как замыкается модель при наличии различных слоев спроса и типов передвижений. После того, как был проведен первый этап модельеру уже известно сколько людей выезжает и въезжает в каждую транспортную зону, а также известно какая доля этих людей совершает поездку того или иного типа, т.е. распределение поездок по слоям спроса. Для каждого слоя спроса формируется своя матрица корреспонденций, при этом используется одна и та же матрица издержек (так как используемая транспортная сеть для всех жителей одна и та же), однако параметры $\{\beta\}$ для каждого слоя спроса свои. После этого, данные матрицы корреспонденций суммируются. Полученная агрегированная матрица корреспонденций используется для проведения этапа расщепления корреспонденций (суммарных) по типам передвижений. Далее алгоритм работает только с ней (агрегированной матрицей) до начала следующей итерации, когда вся процедура повторяется. Более подробно вся процедура рассматривается ниже.

После построения очередной оценки матрицы корреспонденций, используя, опять же, матрицу издержек, происходит расщепление корреспонденций по типам передвижений. Для этого используются модели дискретного выбора, родственные уже описанной выше logit-choice модели. На выходе данного этапа (3-го шага в описанной схеме) получаются матрицы корреспонденций для каждого типа передвижений.

Отметим, что нет единой методологии и стандарта, какую из моделей дискретного выбора стоит использовать. Каждая из моделей имеет свои недостатки. Так, например, для модели logit-choice получаемое стохастическое равновесие даже в простейших постановках может сильно отличаться от равновесия Нэша. Также, результаты данной модели очень чувствительны к тому, как задается набор альтернатив для выбора. Например, расщепление для дерева выбора (личный транспорт, автобус) и для дерева выбора (личный транспорт, зеленый автобус, красный автобус) будут разными (даже если выбор описывается для одной и той же транспортной системы и вся разница только в том, что во втором случае указан цвет автобуса). Наиболее используемыми на практике являются родственные logit-модели Nested Logit Model и Multinomial Logit Model, а также композитная модель Mixed Logit Model (см. [1]).

Важно сказать, что порядок шагов 2 и 3 зависит от того, какой слой спроса описывается, т.е., грубо говоря, от цели поездки и некоторых параметров жителей. Например, слой спроса может определяться как «люди предпенсионного возраста, совершающие поездки из дома на дачу». Иногда слои спроса определяются как цель поездки с учетом пункта отправления, например «поездка из дома на работу». Выше мы отмечали, что характерное время формирования корреспонденций – годы. Это справедливо для поездок дом–работа и обратно. Однако данный тезис кажется не совсем точным, например, для поездок за покупками. В этом случае можно предположить, что люди определяют место покупок уже после того, как будет определен тип передвижения, т.е. для данного слоя спроса, шаги 2 и 3 оказывается возможным поменять местами.

Наконец, последним этапом (шаг 4) является расчет равновесного распределения потоков для каждого типа передвижения. Для личного транспорта при этом используется модель Бэкмана. Для общественного транспорта не существует единого подхода к моделированию и, как правило, используются композитные модели, включающие элементы моделей дискретного выбора. Стоит лишь отметить, что при численном расчете моделей реального города именно этап расчета равновесного распределения потоков является самым затратным по времени. На выходе этого этапа получаются вектора загрузки сети (т.е. потоки по ребрам транспортного графа) и матрица издержек (шаг 5).

Шаги 6–8 в приведенной выше схеме отвечают за создание обратной связи в модели, которая позволит учесть взаимное влияние матрицы корреспонденций и матрицы издержек. Напомним, что мы использовали матрицу издержек для расчета матрицы корреспонденций на шаге 2 и матрицу корреспонденций для получения матрицы издержек на шаге 4. Логично требовать, чтобы матрица издержек, подаваемая на вход на шаге 2 и матрица издержек, получаемая на выходе на шаге 5 если не совпадали, то были «близки» по какому-либо разумному критерию.

Критерием остановки служит равенство (с требуемой точностью) средних издержек для матрицы издержек, полученной по модели расчета матрицы корреспонденций, и средних издержек, полученных эмпирическим путем. Мы к этому еще вернемся чуть позже. В статьях [23, 50] было показано, что данный критерий остановки в нашей постановке соответствует оценке матрицы корреспонденций методом максимального правдоподобия [51] для описанной далее модели.

Пусть имеются данные о реализации некоторой случайной матрицы корреспонденций $\{R_{ij}\}$, каждый элемент которой (независимо от всех остальных) распределен по закону Пуассона с математическим ожиданием $M(R_{ij}) = d_{ij}$, причем верно параметрическое предположение о том, что матрица корреспонденций представляется следующим образом:

$$d_{ij} = \tilde{A}_i \tilde{B}_j f(C_{ij}),$$

где $f(C_{ij})$ – функция притяжения (часто выбирают $f(C_{ij}) = \exp(-\beta C_{ij})$), а матрица издержек $\{C_{ij}\}$ считается известной.[12]

Зададимся целью найти оценку $\{d_{ij}\}$ при заданной (наблюдаемой) матрице $\{R_{ij}\}$ и известной матрице $\{C_{ij}\}$ методом максимума правдоподобия (на основе теоремы Фишера) [51]. Точнее говоря, оценивать требуется не саму матрицу $\{d_{ij}\}$, а неизвестные параметры $\tilde{A}, \tilde{B}, \beta$. Для возможности применять теорему Фишера [51], и таким образом гарантировать хорошие (асимптотические) свойства полученных оценок (в 1-норме), будем считать, что число оцениваемых параметров $p = 2n+1$ и объем выборки $N = \sum_{i,j=1}^{n} R_{ij}$ удовлетворяют следующему соотношению: $p/N \ll 1$.

---

[12] Важно заметить, что мы допускаем, следуя В.Г. Спокойному, “model misspecification” [51], т.е., что эти предположения не верны, и истинное распределение вероятностей $\{R_{ij}\}$ не лежит в этом семействе. Тогда полученное решение по методу максимума правдоподобия можно интерпретировать, как дающее асимптотически (по $N = \sum_{i,j=1}^{n} R_{ij}$ ) наиболее близкое (по расстоянию Кульбака–Лейблера) распределение в этом семействе к истинному распределению. Другими словами, так полученные параметры – являются “асимптотически наилучшими” оценками параметров проекции (по расстоянию Кульбака–Лейблера) истинного распределения вероятностей на выбранное нами параметрическое семейство. Детали см. в [51].

Из определения распределения Пуассона следует, что вероятность реализации корреспонденции $R_{ij}$ может быть посчитана как:

$$P\left(R_{ij} \mid d_{ij}\right) = \frac{\exp\left(-d_{ij}\right) \cdot d_{ij}^{R_{ij}}}{R_{ij}\,!}.$$

Функция правдоподобия (вероятность того, что "выпадет" матрица $\left\{R_{ij}\right\}$, если значения параметров $\tilde{A}, \tilde{B}, \beta$) будет иметь следующий вид:

$$\Lambda\left(\left\{R_{ij}\right\} \middle| \tilde{A}, \tilde{B}, \beta\right) = \prod_{i,j=1}^{n} \frac{\exp\left(-d_{ij}\right) \cdot d_{ij}^{R_{ij}}}{R_{ij}\,!} = \prod_{i,j=1}^{n} \frac{\exp\left(-\tilde{A}_i \tilde{B}_j f\left(C_{ij}\right)\right) \cdot \left(\tilde{A}_i \tilde{B}_j f\left(C_{ij}\right)\right)^{R_{ij}}}{R_{ij}\,!}.$$

Нам нужно найти точку $\left(\tilde{A}, \tilde{B}, \beta\right)$, в которой достигается максимум функции правдоподобия. Как известно, точка максимума неотрицательной функции (каковой по определению является вероятность) не изменится, если решать задачу максимизации не для исходной функции, а для ее логарифма. Перейдем к логарифму функции правдоподобия:

$$\ln \Lambda = \sum_{i,j=1}^{n} \left(-\tilde{A}_i \tilde{B}_j f\left(C_{ij}\right) + R_{ij} \cdot \left(\ln \tilde{A}_i + \ln \tilde{B}_j + \ln f\left(C_{ij}\right)\right) - \ln R_{ij}\,!\right).$$

Получаем следующую задачу оптимизации:

$$\sum_{i,j=1}^{n} \left(-\tilde{A}_i \tilde{B}_j f\left(C_{ij}\right) + R_{ij} \cdot \left(\ln \tilde{A}_i + \ln \tilde{B}_j + \ln f\left(C_{ij}\right)\right) - \ln R_{ij}\,!\right) \to \max_{\tilde{A} \ge 0,\, \tilde{B} \ge 0,\, \beta}.$$

Выпишем условия оптимальности:[13]

$$\frac{\partial \ln \Lambda}{\partial \tilde{A}_i} = \sum_{j=1}^{n} \left(-\tilde{B}_j \exp\left(-\beta C_{ij}\right) + \frac{R_{ij}}{\tilde{A}_i}\right) = 0,$$

$$\frac{\partial \ln \Lambda}{\partial \tilde{B}_j} = \sum_{i=1}^{n} \left(-\tilde{A}_i \exp\left(-\beta C_{ij}\right) + \frac{R_{ij}}{\tilde{B}_j}\right) = 0,$$

$$\frac{\partial \ln \Lambda}{\partial \beta} = \sum_{i,j=1}^{n} \left(\tilde{A}_i \tilde{B}_j C_{ij} \exp\left(-\beta C_{ij}\right) - R_{ij} C_{ij}\right) = 0.$$

Мы получили $2n+1$ уравнение максимума правдоподобия:

$$\sum_{j=1}^{n} \tilde{A}_i \tilde{B}_j \exp\left(-\beta C_{ij}\right) = \sum_{j=1}^{n} R_{ij},$$

[13] Легко понять, что максимум не может достигаться на границе. Если допустить, что, скажем, $L_i$ равно нулю в точке максимума, то поскольку $\partial \ln \Lambda / \partial L_i = \infty$ в этой точке, сдвинувшись немного перпендикулярно гиперплоскости $L_i = 0$ внутрь области определения, мы увеличили бы значения функционала – то есть пришли бы к противоречию с предположением о равенстве нулю $L_i$ в точке максимума.

$$\sum_{i=1}^{n} \tilde{A}_i \tilde{B}_j \exp(-\beta C_{ij}) = \sum_{i=1}^{n} R_{ij},$$

$$\sum_{i,j=1}^{n} \left( \tilde{A}_i \tilde{B}_j C_{ij} \exp(-\beta C_{ij}) - R_{ij} C_{ij} \right) = 0.$$

Положим по определению $L_i = \sum_{j=1}^{n} R_{ij}$ и $W_j = \sum_{i=1}^{n} R_{ij}$, $\tilde{A}_i = L_i A_i$ и $\tilde{B}_j = W_j B_j$. Тогда

$$d_{ij} = A_i L_i B_j W_j f(C_{ij}),$$

где $A_i$, $B_j$ – структурные параметры (гравитационной) модели. Их рассчитывают с помощью простого итерационного алгоритма (метод балансировки = метод простых итераций [6, 18, 19, 29, 30]) по значениям $L_i$, $W_j$ следующим образом:[14]

$$A_i = \frac{1}{\sum_{j=1}^{n} B_j W_j f(C_{ij})}, \quad B_j = \frac{1}{\sum_{i=1}^{n} A_i L_i f(C_{ij})}. \tag{6}$$

Последнее же уравнение системы уравнений максимума правдоподобия даст нам:

$$\sum_{i,j=1}^{n} \left( d_{ij} C_{ij} - R_{ij} C_{ij} \right) = 0 \text{ или } \sum_{i,j=1}^{n} R_{ij} C_{ij} = \sum_{i,j=1}^{n} d_{ij} C_{ij}. \tag{7}$$

**Замечание 2.** Обратим внимание, что к тем же самым соотношениям можно было прийти (при той же параметрической гипотезе), если вместо предположения: $R_{ij}$ – независимые случайные величины, распределенные по законам Пуассона с математическим ожиданием $M(R_{ij}) = d_{ij}$, считать, что $\{R_{ij}\} \gg 1$ имеют мультиномиальное распределение с параметрами $\left\{ d_{ij} \Big/ \sum_{i,j} d_{ij} \right\}$. Такой подход представляется более естественным, чем изложенный выше. Кроме того, поскольку константа строгой выпуклости энтропии в 1-норме равна 1 [13] (неравенство Пинскера, оценивающее снизу расстояние Кульбака–Лейблера с помощью квадрата 1-нормы), то на базе асимптотического представления логарифма функции правдоподобия в окрестности точки максимума [51], можно построить довери-тельный интервал для разности оценок и оцениваемых параметров в (наиболее естественной) 1-норме.

---

[14] При $f(C_{ij}) = \exp(-\beta C_{ij})$ эту модель называют также энтропийной моделью [18] или моделью А.Дж. Вильсона [6, 17]. Энтропийная модель, использованная нами ранее (см. п. 3), для связи матрицы корреспонденций и матрицы издержек является частным случаем гравитационной модели при указанном выборе функции притяжения. Точнее, говоря, если бы мы считали, что $C_{ij} = T_{ij}$ не зависит от $\{d_{ij}\}$ и $R(C_{ij}) = f(C_{ij})/2$, то решение задачи (2) в точности давало бы гравитационную модель. При этом $A_i$, $B_j$ выражались бы через множители Лагранжа (двойственные переменные), соответственно, $\lambda_i^L$, $\lambda_j^W$. Система (6) при этом получалась бы при подстановке решения задачи (2), зависящего от этих $2n$ неизвестных параметров, в ограничения (А) (которых тоже $2n$), точнее говоря, независимых параметров и уравнений было бы $2n$–1 [6].

Соответственно, если критерий остановки (7) не выполняется, то происходит перерасчет матрицы корреспонденций с учетом обновленной матрицы издержек.

Калибровку $\beta$ можно проводить в случае, когда имеется дополнительная информация о реальных издержках на дорогах. Для этого применяется следующий алгоритм.

Пусть $C_l^*$ – средние издержки на проезд в системе (известны, например, из опросов) для $l$-го слоя спроса [1], например, для трудовых миграций.

Алгоритм [52]:

1. Рассчитываем $\beta_l^0 = \frac{1}{C_l^*}$, $m = 0$;
2. Рассчитываем $\left\{d_{ij}\left(\beta_l^m\right)\right\}_{i,j=1}^{n,n}$ – матрицу корреспонденций при $\beta_l = \beta_l^m$;
3. Пересчет «4-стадийной модели»;
4. Рассчитываем $C_l^m$ – средние издержки на проезд, соответствующие матрице корреспонденций $\left\{d_{ij}\left(\beta_l^m\right)\right\}_{i,j=1}^{n,n}$;
5. При $m \geq 1$ проверяем условие: $\left|C_l^{m-1} - C_l^*\right| \leq \varepsilon$ (критерий остановки);
6. Рассчитываем $C_l^1$ по $\beta_l^1 = \frac{\beta_l^0 C_l^0}{C_l^*}$ (при $m = 0$) и полагаем

$$\beta_l^{m+1} = \frac{\left(C_l^* - C_l^{m-1}\right)\beta_l^m + \left(C_l^m - C_l^*\right)\beta_l^{m-1}}{C_l^m - C_l^{m-1}} \text{ (при } m \geq 1);$$

7. Переходим на шаг 2.

Данный алгоритм (и его вариации) дает оценку коэффициента $\beta$ при известных средних издержках в сети и является наиболее часто применяемым на практике [1].

Критерий остановки выбран именно таким по следующей причине. Если критерий остановки из пункта 5) выполняется, то, подставляя $C_l^{m-1} \simeq C_l^*$ в формулу расчета $\beta_l^{m+1}$ на шаге 6, получаем:

$$\beta_l^{m+1} \simeq \frac{\left(C_l^* - C_l^*\right)\beta_l^m + \left(C_l^m - C_l^*\right)\beta_l^{m-1}}{C_l^m - C_l^*} = \beta_l^{m-1}.$$

Полученное значение параметра $\beta\left(C^*\right)$ «соответствует» средним издержкам $C^*$, наблюдаемым «в жизни». Сходимость алгоритма и монотонная зависимость $\beta(C)$ (следовательно, и взаимно-однозначное соответствие) между средними издержками $c$ и значением параметра $\beta$ были показаны в работе [53].

Стоит, однако, отметить, что при численной реализации алгоритма возникает множество проблем. В частности, при «плохом» выборе точки старта или при плохой калибровке параметров функций издержек возникают ситуации, в которых алгоритм требует вычисления экстремально больших чисел и превышает размер доступной памяти даже на самых современных машинах (кластерах).

Другой проблемой является отсутствие оценки скорости сходимости алгоритма. Зачастую, при моделировании критерий остановки устанавливается жестко. Например, прописывается, что алгоритм калибровки $\beta$ должен быть использован не более 20 раз (или другое, разумное, по мнению модельера, количество итераций). При этом отсутствует четкий критерий качества оценок, полученных таким образом.

Еще одной проблемой описанного подхода является высокая чувствительность к параметрам модели с одной стороны, и невозможность учесть реальные подтвержденные данные на отдельных элементах сети с другой. Другими словами, допустим нам известны величины потоков на той или иной дороге каждый день, к сожалению, данную информацию использовать в такой модели не представляется возможным.

## 6. Модель стабильной динамики

Приведем основные положения модели стабильной динамики, следуя [6, 7]. В рамках модели предполагается, что водители действуют оппортунистически, т.е. выполнен первый принцип Вардропа. Рассмотрим ориентированный граф $\Gamma(V,E)$. В модели каждому ребру $e \in E$ ставятся в соответствия параметры $\bar{f}_e$ и $\bar{t}_e$. Они имеют следующую трактовку: $\bar{f}_e$ – максимальная пропускная способность ребра $e$, $\bar{t}_e$ – минимальные временные издержки на прохождение ребра $e$. Таким образом, сама модель задается графом $\Gamma(V,E,\bar{f},\bar{t})$, где $\bar{f} = \left\{\bar{f}_1,...,\bar{f}_{|E|}\right\}^T$, $\bar{t} = \left\{\bar{t}_1,...,\bar{t}_{|E|}\right\}^T$. Пусть $f$ – вектор распределения потоков по ребрам, инициируемый равновесным распределением потоков по маршрутам, а $t$ – вектор временных издержек, соответствующий распределению $f$. Тогда, если транспортная система находится в стабильном состоянии, всегда выполняются неравенства $f \le \bar{f}$ и $t \ge \bar{t}$. При этом считается, что, если поток по ребру $f_e$ меньше, чем максимальная пропускная способность ребра $\bar{f}_e$, то все автомобили в потоке двигаются с максимальной скоростью, а их временные издержки $t_e$ минимальны и равны $\bar{t}_e$. Если же поток по ребру $f_e$ становится равным пропускной способности ребра $\bar{f}_e$, то временные издержки водителей $t_e$ могут быть сколь угодно большими. Это удобно объяснить следующим образом. Допустим, на некоторое ребро $e$ стало поступать больше автомобилей, чем оно способно обслужить. Тогда на этом ребре начинает образовываться очередь (пробка). Временные издержки на прохождение ребра $t_e$ складываются из минимальных временных издержек $\bar{t}_e$ и времени, которое водитель вынужден отстоять в пробке. При этом, очевидно, если входящий поток автомобилей на ребро $e$ не снизится до максимально допустимого уровня (пропускной способности ребра), то очередь будет продолжать расти и система не будет находиться в стабильном состоянии. Если же в какой-то момент входящий на ребро $e$ поток снизится до уровня пропускной способности ребра, то в системе наступит равновесие. При этом пробка на ребре $e$ (если входящий поток $f_e$ будет равен $\bar{f}_e$) не будет рассасываться, т.е. временные издержки так и останутся на уровне $t_e$ ($t_e > \bar{t}_e$). Рассмотрим это на примере из статьи [7].

**Пример 1.** Рассмотрим (в рамках модели стабильной динамики) граф $\Gamma\left(V, E, \bar{t}, \bar{f}\right)$ (см. рис. 2). Пункты 1 и 2 – потокообразующая пара. При этом выполнено: $d_{12}$ – поток из 1 в 2, $\bar{t}_{up} < \bar{t}_{down}$. Если выполнено $d_{12} < \bar{f}_{up}$, то все водители будут использовать ребро up, причем пробка образовываться не будет, так как пропускная способность ребра больше, чем количество желающих проехать (в единицу времени) водителей. В момент, когда $d_{12} = \bar{f}_{up}$ возможности ребра up будут использоваться на пределе. Если же в какой-то момент величина $f$ станет больше $\bar{f}_{up}$, то на ребре *up* начнет образовываться пробка. Пробка будет расти до тех пор, пока издержки от использования маршрута *up* с не сравняются с издержками от использования маршрута *down*. В этот момент оставшаяся часть начнет использовать маршрут *down*. Если же корреспонденция из 1 в 2 превысит суммарную пропускную способность ребер *up* и *down*, то пробки будут расти неограниченно (входящий поток на ребро будет больше, чем исходящий, соответственно количество автомобилей в очереди будет расти постоянно), т.е. стабильное распределение в системе никогда не установится. Более подробно рассмотрение данной задачи (и модели) стоит смотреть в работе [7]. ■

Рассмотрим другой модельный пример из статьи [7].

**Пример 2 (Парадокс Браесса).** Задан граф $\Gamma\left(V, E, \bar{t}, \bar{f}\right)$ (см. рис. 3). При этом выполнено: $\bar{t}_{13} = 1$ час, $\bar{t}_{12} = 15$ минут, $\bar{t}_{23} = 30$ минут; (1,3) и (2,3) – потокообразующие пары, $d_{13} = 1500$ авт/час и $d_{23} = 1500$ авт/час – соответствующие корреспонденции. Пусть пропускные способности всех ребер одинаковы и равны 2000 авт/час. Тогда равновесие будет такое: 500 авт/час из 1 будут направляться в 3 через 2, а 1000 авт/час будут ехать напрямую (для выезжающих из 2 никаких альтернатив нет). При этом все водители, выезжающие из 1, потратят 1 час, а водители, выезжающие из 2, потратят 45 минут. Таким образом, водителям, едущим из 2 в 3 выгодно, чтобы ребро 1–2 отсутствовало, в то время как для водителей, которые едут из 1 в 3, наличие ребра 1–2 безразлично. Т.е., другими словами, если бы мы имели власть запретить проезд по ребру 1–2, то часть водителей выиграла бы от такого решения, и никто бы не проиграл. Интересно заметить, что для модели Бэкмана это не выполняется. Действительно, в модели Бэкмана, как и для модели стабильной динамики, издержки для водителей, следующих из 1 в 3 равны издержкам на ребре 1–3. Однако они монотонно возрастают от потока на данном ребре. Если бы мы запретили проезд по ребру 1–2, то поток на 1–3 увеличился бы, следовательно, возросли бы и издержки на ребре 1–3. Другими словами, улучшение ситуации для водителей, следующих из 2 в 3, привело бы к ухудшению ситуации для водителей, следующих из 1 в 3. ■

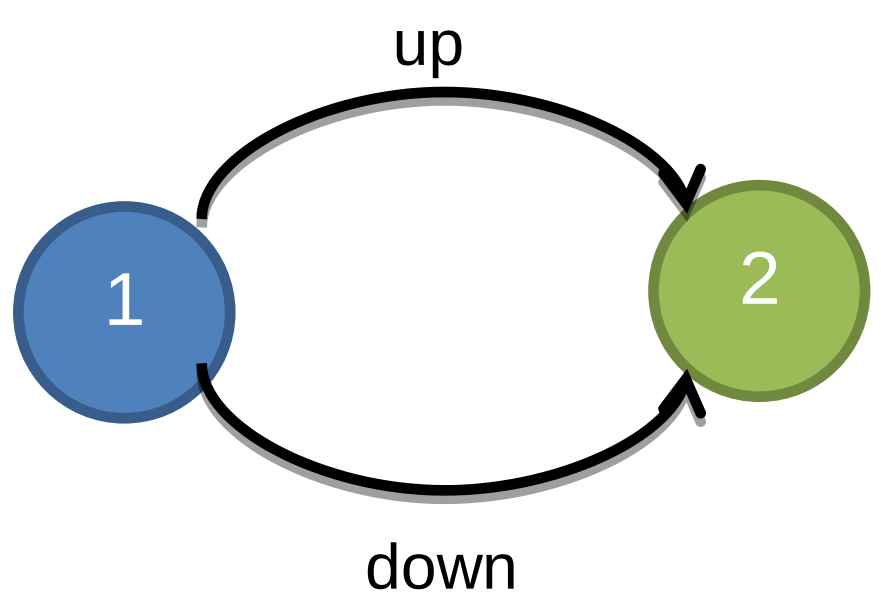

**Рис. 2**

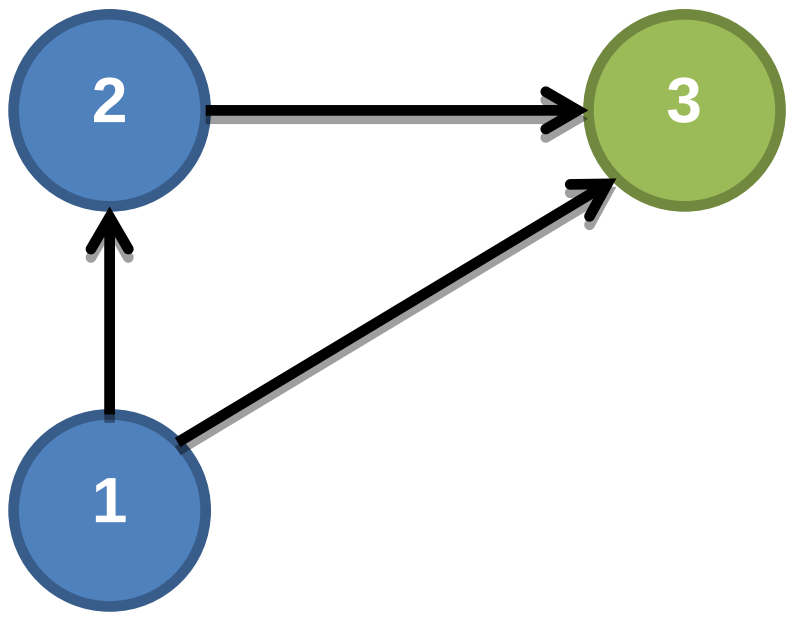

**Рис.3**

Введем ряд новых обозначений. Пусть $w=(i,j)$ – потокообразующая пара (источник – сток) графа $\Gamma(V,E,\bar{t},\bar{f})$, $P_w$ – множество соответствующих $w$ маршрутов, а $t$ – установившийся на графе вектор временных издержек. Тогда временные издержки, соответствующие самому “быстрому” (наикратчайшему) маршруту из $P_w$ равны: $T_w(t)=\min_{p\in P_w}\sum_{e\in E}\delta_{ep}t_e$. Функция $T_w(t)$ зависит от вектора временных издержек $t$. Важно заметить, что функция $T_w(t)$ и её супердифференциал эффективно вычисляются, например, алгоритмом Дейкстры [54, 55] или более быстрыми современными методами [56].

Пусть корреспонденция потокообразующей пары $w$ равна $d_w$. Тогда, если $f^*$ является равновесным распределением потоков для заданного графа $\Gamma(V,E,\bar{t},\bar{f})$, а $t^*$ – соот-ветствующий равновесному распределению вектор временных издержек, то $f^* = \sum_{w\in W} f_w^*$, где $f_w^*$ – вектор распределения потока, порождаемого потокообразующей парой $w$ ($w \in OD$). При этом $f_w^*$ удовлетворяет

$$f_w^* \in d_w \partial T_w(t),$$

где $\partial T_w(t)$ – супердифференциал, который можно посчитать с помощью теоремы Милютина–Дубовицкого [14] используя структуру функции $T_w(t)$. Следовательно, по теореме Моро–Рокафеллара [14] $f^* \in \sum_{w\in W} d_w \partial T_w(t)$. Остается вопрос: как искать равновесный вектор временных издержек $t^*$.

**Теорема 2 [7].** *Распределение потоков $f^*$ и вектор временных издержек $t^*$ являются равновесными для графа $\Gamma(V, E, \bar{t}, \bar{f})$, заданного множества потокообразующих пар $W$ и соответствующих им потоков $d_w$, тогда и только тогда, когда $t^*$ является решением следующий задачи вогнутой оптимизации (максимизация количества свободно движущихся автомобилей, то есть автомобилей не стоящих в пробках):*

$$\max_{t\geq\bar{t}} \left\{ \sum_{w\in W} d_w T_w(t) - \langle \bar{f}, t - \bar{t} \rangle \right\}, \tag{8}$$

*где $f^* = \bar{f} - s^*$, где $s^*$ – (оптимальный) вектор двойственных множителей для ограничений $t \geq \bar{t}$ в задаче (8) или, другими словами, решение двойственной задачи.*

**Схема доказательства.** Резюмируем далее основные постулаты модели стабильной динамики:

1. $t_e \geq \bar{t}_e$, $f_e \leq \bar{f}_e$, $e \in E$;
2. $(t_e - \bar{t}_e)(\bar{f}_e - f_e) = 0$, $e \in E$;
3. $f \in \partial \sum_{w\in W} d_w T_w(t)$ (принцип (Нэша–)Вардропа).

Наша цель, подобрать такую задачу выпуклой (вогнутой) оптимизации, решение которой давало бы пару $(t, f)$, удовлетворяющую этим постулатам. С задачей выпуклой оптимизации намного удобнее работать (разработаны эффективные способы решения таких задач [13]), чем с описанием 1 – 3. Если бы мы могли ограничиться только п. 3, то тогда этот пункт можно было бы переписать

$$\max_{t\geq 0} \left\{ \sum_{w\in W} d_w T_w(t) - \langle f, t - \bar{t} \rangle \right\}, \tag{9}$$

но есть еще пп. 1, 2. Подсказка содержится в п. 2, который удобно понимать как условие дополняющей нежесткости [14] ограничения $t \geq \bar{t}$ (в задаче (9)) с множителями Лагранжа $s = \bar{f} - f \geq 0$. Причем, нет необходимости оговаривать дополнительно, что $f \leq \bar{f}$, поскольку множители Лагранжа к ограничениям вида неравенство автоматически неотрицательны. Таким образом, мы приходим к (8). Все же необходимо оговориться, что получившаяся задача (8) хотя и является довольно простой задачей (не сложно заметить, что это просто задача линейного программирования (ЛП), записанная в более компактной форме, действительно: $\max_{t\geq 0;\{T_w\};C}\left\{C-\left\langle \bar{f},t-\bar{t}\right\rangle;T_w\leq\sum_{e\in E}\delta_{ep}t_e,\ p\in P_w,\ w\in W;C=\sum_{w\in W}d_wT_w\right\}$), но требуется также найти двойственные множители (чтобы определить $f$ ). □

Задача (8) имеет (конечное) решение тогда и только тогда, когда существует хотя бы один вектор потоков $f\leq\bar{f}$ , соответствующий заданным корреспонденциям. Другими словами, существует способ так организовать движение, чтобы при заданных корреспонденциях не нарушались условия $f\leq\bar{f}$ . В противном случае, равновесного (стационарного) режима в системе не будет, и со временем весь граф превратится в одну большую пробку, которая начнет нарушать условия заданных корреспонденций, не позволяя новым пользователям приходить в сеть с той интенсивностью, с которой они этого хотят.

В отличие от модели Бэкмана, в которой при монотонно возрастающих функциях затрат на ребрах равновесие (по потокам на ребрах) единственно, в модели стабильной динамики это может быть не так. Хотя типичной ситуацией (общим положением) будет единственность равновесия, поскольку поиск стабильной конфигурации эквивалентен решению задачи ЛП, все же в определенных “вырожденных” ситуациях может возникать не единственность. Скажем, если в примере 1 (см. рис. 2) $\bar{f}_{up}=d_{12}$. Возникает вопрос: как выбрать единственное равновесие, то которое реализуется? Один из ответов имеется в [7]: следить за динамикой во времени $d_{12}(t)$, а именно за $\int_0^T\left(d_{12}(t)-\bar{f}_{up}\right)dt<\infty$. Другой ответ связан с прямым осуществлением для данного примера предельного перехода, описанного в следующем пункте.

По сути, с помощью этого и следующего пункта мы установим, что если

$$\tau_e^{\mu}(f_e)\xrightarrow[\mu\to 0+]{}\begin{cases}\bar{t}_e, & 0\leq f_e<\bar{f}_e\\ [\bar{t}_e,\infty), & f_e=\bar{f}_e\end{cases},$$

$$d\tau_e^{\mu}(f_e)/df_e\xrightarrow[\mu\to 0+]{}0,\ 0\leq f_e<\bar{f}_e,$$

$x(\mu)$ – равновесное распределение потоков по путям в модели Бэкмана при функциях затрат на ребрах $\tau_e^{\mu}(f_e)$, то

$$\tau_e^{\mu}\left(f_e\left(x(\mu)\right)\right)\xrightarrow[\mu\to 0+]{}t_e,$$

$$f_e\left(x(\mu)\right)\xrightarrow[\mu\to 0+]{}f_e,$$

где пара $(t,f)$ – равновесие в модели стабильной динамики с тем же графом и матрицей корреспонденций, что и в модели Бэкмана, и с ребрами, характеризующимися набором $(\bar{t},\bar{f})$ из определения функций $\tau_e^{\mu}(f_e)$.

Численно решать задачу (1) можно различными способами (см., например, [40]).

## 7. Связь модели стабильной динамики с моделью Бэкмана

Получим модель стабильной динамики в результате предельного перехода из модели Бэкмана. Для этого будем считать, что $\tau_e(f_e)=\bar{t}_e\cdot\left(1-\mu\ln\left(1-f_e/\bar{f}_e\right)\right)$.

Перепишем исходную задачу поиска равновесного распределения потоков (здесь используются обозначения: $\sigma_e^*(t_e)$ – сопряженная функция к $\sigma_e(f_e)$):

$$\min_{f,x}\left\{\sum_{e\in E}\sigma_e\left(f_e\right):\ f=\Theta x,\ x\in X\right\}=\min_{f,x}\left\{\sum_{e\in E}\max_{t_e\in\operatorname{dom}\sigma_e^*}\left[f_e t_e-\sigma_e^*\left(t_e\right)\right]:\ f=\Theta x,\ x\in X\right\}=$$

$$=\max_{t\in\operatorname{dom}\sigma^*}\left\{\min_{f,x}\left[\sum_{e\in E}f_e t_e:\ f=\Theta x,\ x\in X\right]-\sum_{e\in E}\sigma_e^*\left(t_e\right)\right\}.$$

Найдем $\sum_{e\in E}\sigma_e^*\left(t_e\right)$:

$$\sigma_e^*\left(t_e\right)=\sup_{f_e}\left(t_e\cdot f_e-\int_0^{f_e}\tau(z)\,dz\right)=\sup_{f_e}\left(t_e\cdot f_e-\int_0^{f_e}\bar{t}_e\cdot\left(1-\mu\ln\left(1-z/\bar{f}_e\right)\right)dz\right)=$$

$$=\sup_{f_e}\left(\left(t_e-\bar{t}_e\right)f_e+\bar{t}_e\mu\int_0^{f_e}\ln\left(1-z/\bar{f}_e\right)dz\right).$$

Из принципа Ферма [14] найдем $f_e$, доставляющее максимум:

$$\frac{\partial}{\partial f_e}\left(\left(t_e-\bar{t}_e\right)f_e+\bar{t}_e\mu\int_0^{f_e}\ln\left(1-z/\bar{f}_e\right)dz\right)=0\ \Rightarrow$$

$$\exp\left(-\frac{t_e-\bar{t}_e}{\bar{t}_e\mu}\right)=1-f_e/\bar{f}_e\ \Rightarrow$$

$$f_e=\bar{f}_e\cdot\left(1-\exp\left(-\frac{t_e-\bar{t}_e}{\bar{t}_e\mu}\right)\right).$$

Подставляем в $\sigma_e^*\left(t_e\right)$:

$$\left(t_e-\bar{t}_e\right)f_e+\bar{t}_e\mu\int_0^{f_e}\ln\left(1-z/\bar{f}_e\right)dz=\left(t_e-\bar{t}_e\right)f_e-\bar{t}_e\cdot\bar{f}_e\cdot\mu\int_1^{1-f_e/\bar{f}_e}\ln z\,dz=$$

$$=\left(t_e-\bar{t}_e\right)f_e-\bar{t}_e\cdot\bar{f}_e\cdot\mu\left(\left(1-\frac{f_e}{\bar{f}_e}\right)\ln\left(1-\frac{f_e}{\bar{f}_e}\right)-1+\frac{f_e}{\bar{f}_e}+1\right)=$$

$$=\left(t_e-\bar{t}_e\right)f_e-\left(-\left(\bar{f}_e-f_e\right)\left(t_e-\bar{t}_e\right)+f_e\cdot\bar{t}_e\cdot\mu\right)=$$

$$=\left(t_e-\bar{t}_e\right)\bar{f}_e-f_e\cdot\bar{t}_e\cdot\mu=\left(t_e-\bar{t}_e\right)\bar{f}_e-\bar{f}_e\cdot\bar{t}_e\cdot\mu\left(1-\exp\left(-\frac{t_e-\bar{t}_e}{\bar{t}_e\mu}\right)\right).$$

Возвращаясь к исходной задаче имеем:

$$\max_{t\in\operatorname{dom}\sigma^*}\left\{\min_{f,x}\left[\sum_{e\in E}f_e t_e:\ f=\Theta x,\ x\in X\right]-\sum_{e\in E}\sigma_e^*\left(t_e\right)\right\}=$$

$$=\max_{t\in\operatorname{dom}\sigma^*}\left\{\sum_{w\in W}d_w T_w(t)-\left\langle\bar{f},t-\bar{t}\right\rangle+\mu\sum_{e\in E}\bar{f}_e\cdot\bar{t}_e\left(1-\exp\left(-\frac{t_e-\bar{t}_e}{\bar{t}_e\mu}\right)\right)\right\}\overset{\mu\to 0+}{=}$$

$$\overset{\mu\to 0+}{=}\max_{t\geq\bar{t}}\left\{\sum_{w\in W}d_w T_w(t)-\left\langle\bar{f},t-\bar{t}\right\rangle\right\}.$$

Возможность поменять местами порядок взятия максимума и минимума следует из теоремы фон Неймана о минимаксе [15].

Выбор именно функции $\tau_e\left(f_e\right)=\bar{t}_e\cdot\left(1-\mu\ln\left(1-f_e/\bar{f}_e\right)\right)$ не является определяющим. Вместо этой функции может быть взят любой другой гладкий внутренний барьер области $f_e<\bar{f}_e$, например, $\tau_e\left(f_e\right)=\bar{t}_e\cdot\left(1+\mu\dfrac{\bar{f}_e}{\bar{f}_e-f_e}\right)$ – такую функцию вполне естественно использовать при поиске равновесного распределения пользователей сети общественного транспорта [7]. Найдем $\sigma_e^*\left(t_e\right)$:

$$\sigma_e^*\left(t_e\right)=\sup_{f_e}\left(t_e\cdot f_e-\int_0^{f_e}\bar{t}_e\cdot\left(1+\mu\frac{z}{\bar{f}_e-z}\right)dz\right)=\sup_{f_e}\left(\left(t_e-\bar{t}_e\right)\cdot f_e+\bar{t}_e\cdot\bar{f}_e\cdot\mu\cdot\ln\left(1-\frac{f_e}{\bar{f}_e}\right)+\bar{t}_e\cdot f_e\cdot\mu\right).$$

Вновь выписывая условия оптимальности первого порядка имеем:

$$\frac{\partial}{\partial f_e}\left(\left(t_e-\bar{t}_e\right)\cdot f_e+\bar{t}_e\cdot\bar{f}_e\cdot\mu\cdot\ln\left(1-\frac{f_e}{\bar{f}_e}\right)+\bar{t}_e\cdot f_e\cdot\mu\right)=\left(t_e-\bar{t}_e\right)-\bar{t}_e\cdot\mu\cdot\frac{\bar{f}_e}{\bar{f}_e-f_e}+\bar{t}_e\cdot\mu=0\ \Rightarrow$$

$$f_e=\bar{f}_e\cdot\left(1-\frac{\bar{t}_e\cdot\mu}{t_e-\left(1-\mu\right)\bar{t}_e}\right).$$

Подставляя в $\sigma_e^*\left(t_e\right)$, имеем:

$$\sigma_e^*\left(t_e\right)=\left(t_e-\bar{t}_e\right)\cdot\bar{f}_e\cdot\left(1-\frac{\bar{t}_e\cdot\mu}{t_e-\left(1-\mu\right)\bar{t}_e}\right)+\bar{t}_e\cdot\bar{f}_e\cdot\mu\cdot\ln\left(\frac{\bar{t}_e\cdot\mu}{t_e-\left(1-\mu\right)\bar{t}_e}\right)+\bar{t}_e\cdot\mu\cdot\bar{f}_e\cdot\left(1-\frac{\bar{t}_e\cdot\mu}{t_e-\left(1-\mu\right)\bar{t}_e}\right)=$$

$$=\left(t_e-\bar{t}_e\right)\cdot\bar{f}_e+\bar{t}_e\cdot\bar{f}_e\cdot\mu\cdot\ln\left(\frac{\bar{t}_e\cdot\mu}{t_e-\left(1-\mu\right)\bar{t}_e}\right).$$

В итоге получаем:

$$\max_{t\in\operatorname{dom}\sigma^*}\left\{\min_{f,x}\left[\sum_{e\in E}f_e t_e:\ f=\Theta x,\ x\in X\right]-\sum_{e\in E}\sigma_e^*\left(t_e\right)\right\}=$$

$$=\max_{t\in\operatorname{dom}\sigma^*}\left\{\sum_{w\in W}d_w T_w\left(t\right)-\left\langle\bar{f},t-\bar{t}\right\rangle+\mu\sum_{e\in E}\bar{f}_e\cdot\bar{t}_e\cdot\ln\left(1+\frac{t_e-\bar{t}_e}{\bar{t}_e\cdot\mu}\right)\right\}\overset{\mu\to0+}{=}\max_{t\geq\bar{t}}\left\{\sum_{w\in W}d_w T_w\left(t\right)-\left\langle\bar{f},t-\bar{t}\right\rangle\right\}.$$

Наконец, рассмотрим третий вариант выбора функции $\tau_e\left(f_e\right)$:

$$\tau_e\left(f_e\right)=\bar{t}_e\cdot\left(1+\gamma\cdot\left(\frac{f_e}{\bar{f}_e}\right)^{\frac{1}{\mu}}\right).$$

Функции такого вида называются BPR-функциями и наиболее часто применяются при моделировании. Так, при использовании модели Бэкмана обычно полагают $\mu=0.25$, а значение параметра $\gamma$ варьируется от 0.15 до 2 и определяется типом дороги [1]. Найдем $\sigma_e^*\left(t_e\right)$:

$$\sigma_e^*(t_e)=\sup_{f_e}\left(t_e\cdot f_e-\int_0^{f_e}\bar{t}_e\cdot\left(1+\gamma\cdot\left(\frac{z}{\bar{f}_e}\right)^{\frac{1}{\mu}}\right)dz\right)=\sup_{f_e}\left((t_e-\bar{t}_e)\cdot f_e-\bar{t}_e\cdot\gamma\cdot\int_0^{f_e}\left(\frac{z}{\bar{f}_e}\right)^{\frac{1}{\mu}}dz\right)=$$

$$=\sup_{f_e}\left((t_e-\bar{t}_e)\cdot f_e-\bar{t}_e\cdot\frac{\mu}{1+\mu}\cdot\gamma\cdot\frac{f_e^{1+\frac{1}{\mu}}}{\bar{f}_e^{\frac{1}{\mu}}}\right).$$

Аналогично рассуждая:

$$\frac{\partial}{\partial f_e}\left((t_e-\bar{t}_e)\cdot f_e-\bar{t}_e\cdot\frac{\mu}{1+\mu}\cdot\gamma\cdot\frac{f_e^{1+\frac{1}{\mu}}}{\bar{f}_e^{\frac{1}{\mu}}}\right)=t_e-\bar{t}_e-\bar{t}_e\cdot\gamma\cdot\frac{f_e^{\frac{1}{\mu}}}{\bar{f}_e^{\frac{1}{\mu}}}=0\ \Rightarrow\ f_e=\bar{f}_e\cdot\left(\frac{t_e-\bar{t}_e}{\bar{t}_e\cdot\gamma}\right)^{\mu}.$$

Тогда

$$\sigma_e^*(t_e)=\sup_{f_e}\left((t_e-\bar{t}_e)\cdot f_e-\bar{t}_e\cdot\frac{\mu}{1+\mu}\cdot\gamma\cdot\frac{f_e^{1+\frac{1}{\mu}}}{\bar{f}_e^{\frac{1}{\mu}}}\right)=\bar{f}_e\cdot\left(\frac{t_e-\bar{t}_e}{\bar{t}_e\cdot\gamma}\right)^{\mu}\frac{(t_e-\bar{t}_e)}{1+\mu}.$$

В итоге получаем:

$$\max_{t\in\operatorname{dom}\sigma^*}\left\{\sum_{w\in W}d_wT_w(t)-\sum_{e\in E}\bar{f}_e\cdot\left(\frac{t_e-\bar{t}_e}{\bar{t}_e\cdot\gamma}\right)^{\mu}\frac{(t_e-\bar{t}_e)}{1+\mu}\right\}\overset{\mu\to0+}{=}\max_{t\geq\bar{t}}\left\{\sum_{w\in W}d_wT_w(t)-\left\langle\bar{f},t-\bar{t}\right\rangle\right\}.$$

**Замечание 3.** Естественно в контексте всего написанного выше теперь задаться вопросом: а можно ли получить модель стабильной динамики эволюционным образом, то есть подобно тому, как в конце п. 4 была эволюционно проинтерпретирована модель Бэкмана. Действительно, рассмотрим логит динамику (с $T\to0+$) или, скажем, просто имитационную логит динамику [9]. Предположим, что гладкие, возрастающие, выпуклые функции затрат на всех ребрах $\tau_e(f_e)$ при $\mu\to0+$ превращаются в “ступеньки” (многозначные функции)

$$\tau_e(f_e)=\begin{cases}\bar{t}_e, & 0\leq f_e<\bar{f}_e\\ [\bar{t}_e,\infty), & f_e=\bar{f}_e\end{cases}.$$

Согласно п. 4 равновесная конфигурация при таком переходе $\mu\to0+$ должна находиться из решения задачи

$$\sum_{e\in E}\int_0^{f_e}\tau_e(z)dz\to\min_{f=\Theta x,\,x\in X}.$$

Считая, что в равновесии не может быть $\tau_e(f_e)=\infty$ (иначе, равновесие просто не достижимо, и со временем весь граф превратится в одну большую пробку, см. конец п. 6), можно не учитывать в интеграле вклад точек $\bar{f}_e$ (в случае попадания в промежуток интегрирования), то есть переписать задачу следующим образом

$$\min_{f=\Theta x,\,x\in X}\sum_{e\in E}\int_0^{f_e}\left(\bar{t}_e+\delta_{\bar{f}_e}(z)\right)dz\Leftrightarrow\min_{\substack{f=\Theta x,\,x\in X\\ f\leq\bar{f}}}\sum_{e\in E}f_e\bar{t}_e,$$

где

$$\delta_{\bar{f}_e}(z)=\begin{cases}0, & 0\leq z<\bar{f}_e\\ \infty, & z\geq\bar{f}_e\end{cases}.$$

Мы получили задачу линейного программирования. Интересно заметить, что двойственные множители $\eta \geq 0$ к ограничениям $f \leq \bar{f}$ имеют размерность (“физический смысл”) времени (см. интерпретацию двойственных множителей из п. 3). Чему равно это время сразу может быть не ясно из решения только что выписанной задачи линейного программирования. Но если перейти к двойственной задаче, то мы получим уже рассматриваемую нами ранее задачу (9) из п. 6. Причем двойственные множители $\eta$ связаны с временами проезда ребер $t$ следующим образом: $\eta = t - \bar{t} \geq 0$, то есть, действительно, получают содержательную интерпретацию времени, потерянного на ребрах (дополнительно к времени проезда по свободной дороге) из-за наличия пробок. Еще одним подтверждением только что сказанному является то, что если $f_e = \bar{f}_e$, то на ребре $e$ пробка и время прохождения этого ребра $t_e \geq \bar{t}_e$ ($\eta \geq 0$), а если пробки нет, то $f_e < \bar{f}_e$ и $t_e = \bar{t}_e$ ($\eta = 0$). То есть имеют место условия дополняющей нежесткости.

Из данного замечания становится ясно, как оптимально назначать платы за проезд в модели стабильной динамики. Переведем предварительно время в деньги. Назначая платы за проезд по ребрам графа в соответствии с вектором $\eta$, получим, что равновесное распределение пользователей транспортной сети по этой сети будет описываться парой $(f, \bar{t})$, что, очевидно, лучше, чем до введения плат $(f, t)$. Действительно, задачу поиска равновесных потоков $f$:

$$\min_{\substack{f=\Theta x,\, x\in X\\ f\leq \bar{f}}} \sum_{e\in E} f_e \bar{t}_e$$

можно переписать с помощью метода множителей Лагранжа, следующим образом

$$\min_{f=\Theta x,\, x\in X} \sum_{e\in E}\left( f_e \bar{t}_e + \eta_e \cdot \left( f_e - \bar{f}_e \right)\right).$$

Другими словами, у этих двух задач одинаковые решения $f$. Но $t = \bar{t} + \eta$, поэтому тот же самый вектор $f$ будет доставлять решение задаче

$$\min_{f=\Theta x,\, x\in X} \sum_{e\in E} f_e t_e .$$

Поскольку мы знаем, что $f \leq \bar{f}$, то $f$ будет доставлять решение также и этой задаче

$$\min_{\substack{f=\Theta x,\, x\in X\\ f\leq \bar{f}}} \sum_{e\in E} f_e t_e .$$

Таким образом, один и тот же вектор $f$ отвечает оптимальному (с точки зрения социального оптимума) распределению потоков в графе с ребрами, взвешенными согласно векторам $\bar{t}$ и $t$. Более того, не сложно понять, что все это остается верным не только для двух векторов $\bar{t}$ и $t$, но и для целого “отрезка” векторов: $\tilde{t} = \alpha \bar{t} + (1-\alpha) t$, $\alpha \in [0,1]$. К сожалению, если не взимать платы за проезд, то всегда реализуется сценарий $\alpha = 0$, то есть вектор затрат будет $t$ (объяснению причин были посвящены пп. 6 и 7), но если мы взимаем платы согласно вектору $\eta_\alpha = \alpha \cdot (t - \bar{t})$, то вектор реальных затрат (не учитывающих платы) будет $t - \eta_\alpha = (1-\alpha) t + \alpha \bar{t}$. Следовательно, оптимально выбирать $\alpha = 1$, то есть платы согласно $\eta = t - \bar{t}$.

Собственно, всю эту процедуру (назначения плат) можно делать адаптивно: постепенно увеличивая плату за проезд на тех ребрах, на которых наблюдаются пробки. И делать это нужно до тех пор, пока пробки не перестанут появляться.

## 8. Учет расщепления потоков по способам передвижений

Распространим модель стабильной динамики на тот случай, когда (все) пользователи (игроки) транспортной сети имеют возможность выбирать между двумя альтернативными видами транспорта: личным и общественным [49]. Соответственно, теперь у нас имеется информация о $\Gamma_л\left(V,E^л,\bar{t}^л,\bar{f}^л\right)$ – для личного транспорта и $\Gamma_o\left(V,E^o,\bar{t}^o,\bar{f}^o\right)$ – для общественного транспорта. Аналогично рассматривается общий случай. Имеет место следующий результат, получаемый аналогично [7].

**Теорема 3**. *Распределение потоков* $f^*=\left(f^л,f^o\right)$ *и вектор временных издержек* $t^*=\left(t^л,t^o\right)$ *являются равновесными для графа* $\Gamma\left(V,E,\bar{t},\bar{f}\right)$*, заданного множества потокообразующих пар* $W$ *и соответствующих им потоков* $d_w$*, тогда и только тогда, когда* $t^*$ *является решением следующий задачи выпуклой оптимизации (в дальнейшем для нас будет удобно писать эту задачу как задачу минимизации, а не максимизации, как раньше):*

$$\min_{\substack{t^л\ge\bar{t}^л\\ t^o\ge\bar{t}^o}}\left\{\left\langle\bar{f}^л,t^л-\bar{t}^л\right\rangle+\left\langle\bar{f}^o,t^o-\bar{t}^o\right\rangle-\sum_{w\in W}d_w\min\left\{T_w^л\left(t^л\right),T_w^o\left(t^o\right)\right\}\right\},\qquad(10)$$

*где* $f^л=\bar{f}^л-s^л$, $f^o=\bar{f}^o-s^o$, *а* $s^л$, $s^o$ *– (оптимальные) векторы двойственных множителей для ограничений* $t^л\ge\bar{t}^л$, $t^o\ge\bar{t}^o$ *в задаче (10) или, другими словами, решение двойственной задачи.*

Если мы хотим учитывать возможность пересаживания пользователей сети в пути с личного транспорта на общественный (и наоборот), то выражение $\min\left\{T_w^л\left(t^л\right),T_w^o\left(t^o\right)\right\}$ нужно будет немного изменить. Мы не будем здесь вдаваться в детали, скажем лишь, что с точки зрении всего дальнейшего это не принципиально. Более того, хотя на время работы алгоритма это и скажется (время увеличится), но не критическим образом, поскольку решение задачи о кратчайшем маршруте, которое возникает на каждом шаге субградиентного спуска, естественным образом (небольшим раздутием графа) обобщается на случай когда есть несколько весов ребер, отвечающих разным типам транспортных средств, и в вершинах графа есть затраты на пересадку (изменение транспортного средства), при условии наличия возможности её осуществления.

Отметим, что единого подхода для моделирования общественного транспорта пока не существует. Указанный нами метод соответствует подходу, когда предполагается, что пассажиры выбирают только “оптимальные стратегии”, т.е. соответствует концепции равновесия по Нэшу. Другой, альтернативный подход, который сейчас чаще используется на практике, следует концепции «стохастического равновесия». В нем предполагается, что пассажиры выбирают каждый из возможных маршрутов с некоторой вероятностью, зависящей от ожидаемых издержек в сети и ожидаемого времени ожидания соответствующего маршрута на остановочной станции. Об этом немного подробнее будет рассказано в п. 10.

Обратим внимание на то, что общественный транспорт все же более естественно описывать моделью Бэкмана, а не моделью стабильной динамики. Во всяком случае, к та-кому выводу склоняются авторы модели стабильной динамики [7]. В этой связи далее приводится смешанная модель, в которой личный транспорт описывается моделью стабильной динамики, а общественный – моделью Бэкмана (см. п. 7):

$$\min_{\substack{t^л \ge \bar{t}^л \\ t^o \ge \bar{t}^o\cdot(1-\mu)}} \left\{ \left\langle \bar{f}^л, t^л - \bar{t}^л \right\rangle + \left\langle \bar{f}^o, t^o - \bar{t}^o \right\rangle - \mu \sum_{e\in E} \bar{f}_e^o \cdot \bar{t}_e^o \cdot \ln\left(1 + \frac{t_e^o - \bar{t}_e^o}{\bar{t}_e^o \cdot \mu}\right) - \sum_{w\in W} d_w \min\left\{T_w^л\left(t^л\right), T_w^o\left(t^o\right)\right\} \right\},$$

где $f^л = \bar{f}^л - s^л$, $f_e^o = \bar{f}_e^o \cdot \left(1 - \frac{\bar{t}_e^o \cdot \mu}{t_e^o - (1-\mu)\bar{t}_e^o}\right)$, $s^л$ – (оптимальный) вектор двойственных множителей для ограничений $t^л \ge \bar{t}^л$.

Несложно распространить все проводимые в следующих пунктах рассуждения именно на такую версию модели стабильной динамики с учетом расщепления по типам передвижений. Однако мы не будем здесь этого делать. Укажем лишь, что для разработки эффективного способа решения возникающих задач оптимизации потребуется привлекать прямо-двойственные субградиентные спуски для задач композитной оптимизации [57].

# 9. Трехстадийная модель стабильной динамики

Объединим теперь задачи (2) и (10) в одну задачу:[15]

$$\min_{\lambda^L,\lambda^W} \max_{\sum_{i,j=1}^n d_{ij}=1,\, d_{ij}\ge 0} \left[ -\sum_{i,j=1}^n d_{ij}\ln d_{ij} + \sum_{i=1}^n \lambda_i^L \left(l_i - \sum_{j=1}^n d_{ij}\right) + \sum_{j=1}^n \lambda_j^W \left(\mathrm{w}_j - \sum_{i=1}^n d_{ij}\right) + \right.$$

$$\left. + \beta \min_{\substack{t^л \ge \bar{t}^л \\ t^o \ge \bar{t}^o}} \left\{ \left\langle \bar{f}^л, t^л - \bar{t}^л \right\rangle + \left\langle \bar{f}^o, t^o - \bar{t}^o \right\rangle - \sum_{i,j=1}^n d_{ij} \min\left\{T_{ij}^л\left(t^л\right), T_{ij}^o\left(t^o\right)\right\} \right\} \right] =$$

$$= \min_{\substack{t^л \ge \bar{t}^л \\ t^o \ge \bar{t}^o \\ \lambda^L,\lambda^W}} \left\{ \max_{\sum_{i,j=1}^n d_{ij}=1,\, d_{ij}\ge 0} \left[ -\sum_{i,j=1}^n d_{ij}\ln d_{ij} - \beta \sum_{i,j=1}^n d_{ij} \min\left\{T_{ij}^л\left(t^л\right), T_{ij}^o\left(t^o\right)\right\} - \right. \right.$$

$$\left. \left. - \sum_{i,j=1}^n d_{ij}\cdot\left(\lambda_i^L + \lambda_j^W\right) \right] + \sum_{i=1}^n \lambda_i^L l_i + \sum_{j=1}^n \lambda_j^W \mathrm{w}_j + \beta\left\langle \bar{f}^л, t^л - \bar{t}^л \right\rangle + \beta\left\langle \bar{f}^o, t^o - \bar{t}^o \right\rangle \right\} =$$

$$= \min_{\substack{t^л \ge \bar{t}^л \\ t^o \ge \bar{t}^o \\ \lambda^L,\lambda^W}} \left\{ \ln\left( \sum_{i,j=1}^n \exp\left(-\beta \min\left\{T_{ij}^л\left(t^л\right), T_{ij}^o\left(t^o\right)\right\} - \lambda_i^L - \lambda_j^W\right) \right) + \right.$$

$$\left. + \sum_{i=1}^n \lambda_i^L l_i + \sum_{j=1}^n \lambda_j^W \mathrm{w}_j + \beta\left\langle \bar{f}^л, t^л - \bar{t}^л \right\rangle + \beta\left\langle \bar{f}^o, t^o - \bar{t}^o \right\rangle \right\}, \quad (11)$$

причем $d_{ij} = \tilde{Z}^{-1}\exp\left(-\beta\min\left\{T_{ij}^л\left(t^л\right), T_{ij}^o\left(t^o\right)\right\} - \lambda_i^L - \lambda_j^W\right)$, где $\tilde{Z}^{-1}$ – ищется из условия нормировки, $f^* = \left(f^л, f^o\right)$, $f^л = \bar{f}^л - s^л$, $f^o = \bar{f}^o - s^o$, а $s^л$, $s^o$ – (оптимальные) векторы двойственных множителей для ограничений $t^{л} \ge \bar{t}^{л}$, $t^{o} \ge \bar{t}^{o}$ в задаче (11). Таким образом, все что осталось сделать, это решить задачу негладкой выпуклой оптимизации (11) прямо-двойственным методом. Заметим при этом, что на переменные $t^{л}$, $t^{o}$ – ограничения сверху возникают из вполне естественных соображений. Пусть $\bar{f}$ максимально возможный поток на ребре, поскольку нас интересует оценка сверху, то можно считать $t > \bar{t}$, стало быть, $f = \bar{f}$. Представим себе самую плохую ситуацию: все ребро стоит в пробке. Пусть длина ребра $L$, а средняя длина автомобиля $l$, число полос $r$. Тогда $t \le Lr\big/\left(l\bar{f}\right) + \bar{t}$ не может превышать нескольких часов. Задачу (11) можно решать, например, прямо-двойственным универсальным композитным методом треугольника [58] или методом из работ [59, 60].

[15] Отметим, что все приводимые далее выкладки (а также выкладки следующего пункта) можно провести, отталкиваясь не от задачи (2), а от её упрощенного варианта, описанного в сноске 9.

Отметим, что если известна информация о потоках по ряду дуг $f_e^{л} = \tilde{f}_e^{л}$, $e \in \mathrm{A}^{л}$; $f_e^{o} = \tilde{f}_e^{o}$, $e \in \mathrm{A}^{o}$, то эту информацию можно "зашить" в модель (3С) подобно тому, как это делается в работе [7]: а именно брать минимум по множеству $t_e^{л} \ge \bar{t}_e^{л}$, $e \in E \setminus \mathrm{A}^{л}$, $t_e^{л} \ge 0$, $e \in \mathrm{A}^{л}$, $t_e^{o} \ge \bar{t}_e^{o}$, $e \in E \setminus \mathrm{A}^{o}$, $t_e^{o} \ge 0$, $e \in \mathrm{A}^{o}$, а слагаемые $+\beta\left\langle \bar{f}^{л}, t^{л} - \bar{t}^{л}\right\rangle + \beta\left\langle \bar{f}^{o}, t^{o} - \bar{t}^{o}\right\rangle$ стоит заменить на (параметры $\bar{t}_e^{л}$, $e \in \mathrm{A}^{л}$ и $t_e^{o}$, $e \in \mathrm{A}^{o}$ – неизвестны, но они и не нужны для расчетов, поскольку входят в виде аддитивных констант в функционал)

$$+\beta \sum_{e \in E \setminus \mathrm{A}^{л}} \bar{f}_e^{л} \cdot \left(t_e^{л} - \bar{t}_e^{л}\right) + \beta \sum_{e \in \mathrm{A}^{л}} \tilde{f}_e^{л} \cdot \left(t_e^{л} - \bar{t}_e^{л}\right) + \beta \sum_{e \in E \setminus \mathrm{A}^{o}} \bar{f}_e^{o} \cdot \left(t_e^{o} - \bar{t}_e^{o}\right) + \beta \sum_{e \in \mathrm{A}^{o}} \tilde{f}_e^{o} \cdot \left(t_e^{o} - \bar{t}_e^{o}\right).$$

# 10. Стохастический вариант трехстадийной модели стабильной динамики

Рассуждая аналогично тому, как мы делали выше, можно обобщить результаты п. 9 на случай, когда вместо модели стабильной динамики используется её стохастический вариант (с параметром $T > 0$) [6, 38, 61]:

$$\min_{\substack{t^{л} \ge \bar{t}^{л} \\ t^{o} \ge \bar{t}^{o} \\ \lambda^{L}, \lambda^{W}}} \left\{ \ln\left( \sum_{i,j=1}^{n} \exp\left( \beta T \psi_{ij}\left( \frac{t^{л}}{T}, \frac{t^{o}}{T} \right) - \lambda_i^{L} - \lambda_j^{W} \right) \right) + \right.$$

$$\left. + \sum_{i=1}^{n} \lambda_i^{L} l_i + \sum_{j=1}^{n} \lambda_j^{W} \mathrm{w}_j + \beta\left\langle \bar{f}^{л}, t^{л} - \bar{t}^{л} \right\rangle + \beta\left\langle \bar{f}^{o}, t^{o} - \bar{t}^{o} \right\rangle \right\}, \qquad (12)$$

где

$$\psi_{ij}\left(t^{л}, t^{o}\right) = \ln\left( \sum_{p \in P_{(i,j)}^{Л}} \exp\left( -\sum_{e \in E} \delta_{ep} t_e^{л} \right) + \sum_{p \in P_{(i,j)}^{O}} \exp\left( -\sum_{e \in E} \delta_{ep} t_e^{o} \right) \right).$$

Причем здесь, также как и в предыдущем пункте,

$$d_{ij} = \breve{Z}^{-1} \exp\left(\beta T \psi_{ij}\left(t^{л}/T, t^{o}/T\right) - \lambda_i^L - \lambda_j^W\right),$$

где $\breve{Z}^{-1}$ – ищется из условия нормировки, $f^* = \left(f^{л}, f^{o}\right)$, $f^{л} = \bar{f}^{л} - s^{л}$, $f^{o} = \bar{f}^{o} - s^{o}$, а $s^{л}$, $s^{o}$ – (оптимальные) векторы двойственных множителей для ограничений $t^{л} \geq \bar{t}^{л}$, $t^{o} \geq \bar{t}^{o}$ в задаче (11). Таким образом, все что осталось сделать, это решить задачу гладкой выпуклой оптимизации (12) прямо-двойственным методом, например, [58–60]. Мы не будем здесь подробно описывать возможные способы решения. Обратим только внимание на то, что задача вычисления значений и градиента функции $\psi\left(t^{л}, t^{o}\right)$ – вычислительно сложная, и, на первый взгляд, даже кажется бесперспективной из-за потенциально экспоненциально большого числа возможных маршрутов. Однако с помощью "аппарата характеристических функций на ориентированных графах" [6, 61, 44] и быстрого автоматического дифференцирования [62, 63] пересчитывать значения функции $\psi\left(t^{л}, t^{o}\right)$ и ее градиента можно довольно эффективно – см. алгоритм из работ [6, 61, 44], который вырождается при $T \to 0+$ в алгоритм Беллмана–Форда [54–56]. Отметим, что при этом предельном переходе модель п. 10 перейдет в модель п. 9.

Отметим также, что все сказанное здесь переносится и на функции $\psi_{ij}\left(t^{л}, t^{o}\right)$ более общего вида, соответствующие различным иерархическим способам выбора (типа Nested Logit) [38, 41]. Например, "практически бесплатно" можно сделать такую замену:

$$T\psi_{ij}\left(\frac{t^{л}}{T}, \frac{t^{o}}{T}\right) \to$$

$$\to \eta \ln\left(\exp\left(T \ln\left(\sum_{p \in P_{(i,j)}^{Л}} \exp\left(-\sum_{e \in E} \delta_{ep} \frac{t_e^{л}}{T}\right)\right) \Big/ \eta\right) + \exp\left(T \ln\left(\sum_{p \in P_{(i,j)}^{O}} \exp\left(-\sum_{e \in E} \delta_{ep} \frac{t_e^{o}}{T}\right)\right) \Big/ \eta\right)\right).$$

## 11. Калибровка модели стабильной динамики

О практическом использовании модели стабильной динамики написано, например, в [64] и немного в [7]. В этом пункте мы сконцентрируем внимание на потенциальных плюсах описанной модели с точки зрения её калибровки. Далее мы опускаем нижней индекс $e$.

Способ оценки $\bar{t}$ довольно очевиден: $\bar{t}$ определяется длиной участка (ребра) и типом ребра (в сельской местности, в городе, шоссе). Вся эта информация обычно бывает доступной. С оценкой $\bar{f}$ немного посложнее. Пусть в конце дуги, состоящей из $r$ полос, стоит светофор, который пускает поток с этой дуги долю времени $\chi$. Пусть $q_{\max} \approx 1800\left[авт/час\right]$ – максимально возможное значение потока по одной полосе (это понятие не совсем корректное, но в первом приближении, им можно пользоваться [6], и оно довольно универсально). Тогда [6] $\bar{f} \approx \chi r q_{\max}$. Тут имеется важный нюанс. Во время зеленой фазы (в зависимости от того, что это за фаза, скажем, движение прямо и налево или движение только направо) “работают” не все полосы: $\bar{f} \approx \chi_1 r_1 q_{\max} + \ldots + \chi_l r_l q_{\max}$, где $l$ – число фаз, “пускающих” поток с рассматриваемой дуги, $\chi_k$ – доля времени отводимое фазе $k$, $r_k$ – эффективное число полос рассматриваемой дуги, задействованных на фазе $k$. Раздобыть информацию по работе светофоров как правило не удается. Однако, существует определенные регламенты, согласно которым и устанавливаются фазы работы светофора. Поскольку все довольно типизировано, то часто бывает достаточно иметь информацию только о полосности дорог. Такая информация с 2012 года уже включается в коммерческие системы ГИС.

К сожалению, в реальных транспортных сетях в часы пик бывают пробки, которые приводят к тому, что пропускная способность ребра определяется пропускной способностью не данного ребра, а какого-то из впереди идущих (по ходу движения) ребер. Другими словами, пробка с ребра полностью “заполнила” это ребро и распространилась на ребра, входящие в это ребро. В таких ситуациях $\bar{f}$ нужно считать исходя из ограничений на пропускную способность на впереди идущих ребрах. При этом в выборе разбиения транспортного графа на ребра и в определении на этих ребрах значений $\bar{f}$, $\bar{t}$ стоит исходить из того, чтобы в типичной ситуации пробка если и переходила с ребра на ребро, то, желательно, чтобы это происходило без изменения пропускной способности входящего ребра. Этого не всегда можно добиться, поскольку часто приходится осуществлять разбиение без особого произвола, исходя из въездов/съездов, перекрестков. В таких случаях, формула $\bar{f} \approx \chi_1 r_1 q_{\max} + \ldots + \chi_l r_l q_{\max}$ является лишь оценкой сверху для “среднего” (типичного) значения, которое надо подставлять в модель (часть этого бремени придется переносить и на $\bar{t}$, увеличивая его).

Заметим, что также как и для обычных многостадийных моделей, для калибровки предложенной модели необходимо использовать один и тот же промежуток времени (скажем, обеденные часы) каждого типового дня (например, буднего) – в зависимости от целей. Кроме того, для Москвы в утреннее часы пик значения потоков существенно не стационарные, и в течение час они, как правило, меняются сильно. Вместе с этим среднее время в пути оказывается больше часа. Таким образом, возникает вопрос: что понимается, например, под равновесным распределением потоков $f$? Ответ: “средние значения” нуждается в пояснении. Средние значения не только по дням, но и по исследуемому промежутку времени (в несколько часов) внутри каждого дня. Как следствие получаем, что выводы, сделанные по предложенной модели, нельзя, например, напрямую использовать для краткосрочного прогнозирования ситуации на дорогах или адаптивного управления светофорной сигнализацией. Таким образом, модель работает и выдает не реальные данные, хотя многие переменные модели и обозначают реальные физические параметры транспортного потока, а лишь некоторые средние (агрегированные) показатели, которые, тем не менее, многое могут сказать о ситуации на дорогах.

В данной работе была предложена модель, описывающая равновесие в транспортной сети с точки зрения макро масштабов времени (месяцы). Для исследования поведения транспортного потока внутри одного дня требуются микро модели. Контекст, в котором такие модели используются, часто связан с необходимостью краткосрочного прогнозирования (на несколько часов вперед) и оптимального управления, например, светофорами. Как правило, в таких микро моделях численно исследуется начально-краевая задача для нелинейных УЧП (системы законов сохранения), в которой начальные условия имеются в наличии, а краевые условия (характеристики источников и стоков во времени, матрицы перемешивания в узлах графа транспортной сети) определяются исходя из исторической информации. Однако в последнее время в Москве все чаще можно слышать предложения о том, чтобы максимально информировать участников дорожного движения с целью лучшей маршрутизации. Такая полная информированность приводит к необходимости определять, скажем, матрицы перемешивание не из (и не только из) исторической информации, но и исходя из равновесных принципов, заложенных в модели Бэкмана и стабильной динамики. Причем нам представляется, что модель стабильной динамики подходит на много лучше для этих целей (мотивация имеется в работе [7], в которой на простом примере показывается, как возникает равновесная конфигурация модели стабильной динамики из внутредневной динамики водителей). Как отмечалось в работе [7], за исключением модельных примеров (рис. 2), очень сложной с вычислительной точки зрения представляется задача описания динамического режима функционирования модели, имеющей в своей основе модель стабильной динамики (или какую-то другую равновесную модель). Более того, пока, насколько нам известно, не было предложено ни одной более менее обоснованной модели такого типа (иногда такие модели называют моделями динамического равновесия [1]). Здесь мы лишь укажем на одно, на наш взгляд, перспективное направление: объединить модель стабильной динамики с микромоделью CTM К. Даганзо [6] в современном ее варианте [65].

# СПИСОК ЦИТИРОВАННОЙ ЛИТЕРАТУРЫ